\documentclass[11pt]{article}
\usepackage{amsmath,amssymb}
\usepackage{amsmath,amsfonts,amssymb}
\usepackage{cancel}
\usepackage{shadow}
\usepackage{fancybox}
\usepackage{color}
\usepackage[dvipsnames]{xcolor}
\usepackage{stmaryrd}
\usepackage{dsfont}

\usepackage{anysize} 
\usepackage{indentfirst} 
\usepackage{fancyhdr}

\usepackage{enumerate} 
\newcommand\Aflat{\AA^\flat}
\renewcommand{\Pr}{\operatorname{\textsc{\textbf{Pr}}}}

\newcommand{\Alg}{\operatorname{\textsc{Alg}}}

\newcommand{\BrTens}{\operatorname{\textsc{\textbf{BrTens}}}}
\newcommand{\Tens}{\operatorname{\textsc{\textbf{Tens}}}}
\newcommand{\Bord}{\operatorname{\textsc{\textbf{Bord}}}}

\newcommand{\DGGPR}{\operatorname{DGGPR}}
\newcommand{\WRT}{\operatorname{WRT}}

\DeclareMathOperator{\Triv}{Triv}
\DeclareMathOperator{\Hollow}{Hollow}

\DeclareMathOperator{\Cat}{Cat}

\DeclareMathOperator{\Bimod}{Bimod}

\newcommand{\Vect}{\mathrm{Vect}}
\newcommand{\Proj}{\operatorname{Proj}}

\newcommand{\cob}{\operatorname{\textbf{Cob}}}

\newcommand{\unit}{\ensuremath{\mathds{1}}}
\newcommand{\Id}{\operatorname{Id}}

\newcommand{\FK}{{\Bbbk}}
\newcommand\simdual{\doteq}

\newcommand\surj{\twoheadrightarrow}

\newcommand\Dd{\mathbb{D}}

\newcommand\Cc{\mathbb{C}}
\newcommand\Nn{\mathbb{N}}

\renewcommand\AA{\mathcal{A}}

\newcommand\BB{\mathcal{B}}
\newcommand\CC{\mathcal{C}}
\newcommand\DD{\mathcal{D}}
\newcommand\EE{\mathcal{E}}
\newcommand\FF{\mathcal{F}}

\newcommand\MM{\mathcal{M}}
\newcommand\NN{\mathcal{N}}

\newcommand\PP{\mathcal{P}}

\newcommand\RR{\mathcal{R}}
\newcommand\ZZ{\mathcal{Z}}

\newcommand\VV{\mathcal{V}}
\newcommand\Vv{\mathscr{V}}
\newcommand\WW{\mathcal{W}}

\newcommand{\rcup}{\tikz[scale = 0.3] \draw[->] (0,0) arc(-180:0:1);}
\newcommand{\lcap}{\tikz[scale = 0.3] \draw[<-] (0,0) arc(180:0:1);}

\newcommand{\mt}{{\operatorname{\mathsf{t}}}}

\usepackage[standard, thref, framed, thmmarks]{ntheorem}
\usepackage{framed}

\theoremsymbol{$\lrcorner$}
\theorembodyfont{}
\theorempreskip{6pt}
\theorempostskip{6pt}
\theoremseparator{:}
\newtheorem{mydef}{$\ulcorner$ Definition}[section]

\theoremsymbol{}

\theorembodyfont{\itshape}
\theoremseparator{:}
\newshadedtheorem{myth}[mydef]{Theorem}
\newshadedtheorem{mythq}[mydef]{Conjecture}
\newshadedtheorem{conjecture}[mydef]{Conjecture}
\newshadedtheorem{hypothesis}[mydef]{Hypothesis}

\theoremsymbol{}

\theorembodyfont{\itshape}
\theoremseparator{:}
\newshadedtheorem{prop}[mydef]{Proposition}
\newshadedtheorem{mylemma}[mydef]{Lemma}
\newshadedtheorem{corr}[mydef]{Corollary}

\theorembodyfont{}
\newtheorem{myex}[mydef]{Example }

\theoremheaderfont{}
\newtheorem{rmk}[mydef]{\textit{Remark }}
\theoremsymbol{$\square$}
\newtheorem*{myproof}{\sc Proof }

\theoremstyle{margin}
\theoremsymbol{}
\theoremseparator{}

\usepackage{tikz}
\usetikzlibrary{arrows,decorations,positioning}
\usepackage{tikz-cd}
\usepackage{wrapfig}
\usepackage[top=2cm, bottom=3cm, left=2.5cm, right=2.5cm]{geometry}
\usepackage{setspace}
\usepackage[toc,page]{appendix}
\usepackage{titlesec}
\usepackage{mathrsfs}
\usepackage{cite}
\usepackage{graphicx}
\usepackage{bm}
\usepackage{calc}
\usepackage{hyperref} 
\hypersetup{pdfstartview=XYZ}

\newcommand{\Ind}{\operatorname{Ind}}
\newcommand{\Free}{\operatorname{Free}}
\newcommand{\colim}{\operatorname{colim}}
\newcommand{\Hom}{\operatorname{Hom}}
\newcommand{\Fun}{\operatorname{Fun}}
\newcommand{\End}{\operatorname{End}}

\begin{document}

\newcommand{\daj}[1]{\textcolor{red}{DJ: #1}}
\newcommand{\FC}[1]{\textcolor{red}{FC: #1}}
\newcommand{\BH}[1]{\textcolor{blue}{BH: #1}}
\newcommand{\VVrpz}{\begin{tikzpicture}[baseline = 10pt]
\fill[blue!30] (0,0) rectangle (1,1);
\end{tikzpicture}}
\newcommand{\AArpz}{\begin{tikzpicture}[baseline = -20pt, yscale=-1]
\fill[blue!30] (0,0.5) rectangle (1,1);
\fill[gray!20] (0,0) rectangle (1,0.5);
\draw[very thick] (0,0.5)--(1,0.5);
\end{tikzpicture}}
\newcommand{\AAoprpz}{\begin{tikzpicture}[baseline = 10pt]
\fill[blue!30] (0,0.5) rectangle (1,1);
\draw[very thick] (0,0.5)--(1,0.5);
\end{tikzpicture}}
\newcommand{\CCrpz}{\begin{tikzpicture}[baseline = 10pt]
\fill[blue!30] (0,0) rectangle (1,0.5); 
\draw[very thick] (0,0.5)--(1,0.5); 
\fill[blue!30] (0,0.75) rectangle (1,1.25); 
\draw[very thick] (0,0.75)--(1,0.75);\end{tikzpicture}}
\newcommand{\DDrpz}{\begin{tikzpicture}[baseline = 10pt] 
\fill[blue!30] (0,0) rectangle (1,1); 
\draw[very thick] (0,0)--(1,0); \draw[very thick] (0,1)--(1,1);\end{tikzpicture}}
\newcommand{\IdVrpz}{\begin{tikzpicture}[baseline = 10pt]
\fill[blue!30] (0,0) rectangle (1,1.25);
\end{tikzpicture}}
\newcommand{\Mrpz}{\begin{tikzpicture}[baseline = 10pt] 
\fill[blue!30](0,0.5)--(0,0)--(2,0)--(2,1.25)--(0,1.25)--(0,0.75)--(1,0.75)--(1,0.5)--cycle;
\draw[very thick] (0,0.75)--(1,0.75);
\draw[very thick] (0,0.5)--(1,0.5);
\end{tikzpicture}}
\newcommand{\Mbarrpz}{\begin{tikzpicture}[baseline = 10pt, xscale = -1] 
\fill[blue!30](0,0.5)--(0,0)--(2,0)--(2,1.25)--(0,1.25)--(0,0.75)--(1,0.75)--(1,0.5)--cycle;
\draw[very thick] (0,0.75)--(1,0.75);
\draw[very thick] (0,0.5)--(1,0.5);
\end{tikzpicture}}
\newcommand{\Nrpz}{\begin{tikzpicture}[baseline = 10pt] 
\fill[blue!30] (0.5,0) rectangle (1,1);
\draw[very thick] (1,1)--(0.5,1)--(0.5,0)--(1,0);
\end{tikzpicture}}
\newcommand{\Nbarrpz}{\begin{tikzpicture}[baseline = 10pt, xscale = -1] 
\fill[blue!30] (0.5,0) rectangle (1,1);
\draw[very thick] (1,1)--(0.5,1)--(0.5,0)--(1,0);
\end{tikzpicture}}

\title{Unit inclusion in a non-semisimple braided tensor category and non-compact relative TQFTs}
\author{Benjamin Ha\"ioun }
\date{\vskip-5pt}
\maketitle
\begin{abstract}
The inclusion of the unit in a braided tensor category $\VV$ induces a 1-morphism in the Morita 4-category of braided tensor categories
$\BrTens$. We give criteria for the dualizability of this morphism.

When $\VV$ is a semisimple (resp. non-semisimple) modular category, we show that the unit inclusion induces under the Cobordism Hypothesis a (resp. non-compact) relative 3-dimensional topological quantum field theory. Following Jordan--Safronov, we conjecture that these relative field theories together with their bulk theories recover Witten--Reshetikhin--Turaev (resp. De Renzi--Gainutdinov--Geer--Patureau-Mirand--Runkel) theories, in a fully extended setting. In particular, we argue that these theories can be obtained by the Cobordism Hypothesis.
\end{abstract} 
\setcounter{tocdepth}{2}
\tableofcontents

\section{Introduction}
This paper is motivated by the quest to bridge topological and higher-categorical constructions of Topological Quantum Field Theories (TQFTs). In the first approach one explicitly defines an $n$-manifold invariant and works their way to a TQFT, adding structure or extra conditions as necessary. This is the approach behind Reshetikhin--Turaev's construction \cite{TuraevBook} of the 3-TQFTs predicted by Witten \cite{WittenJonesPol}, and their non-semisimple variants \cite{DGGPR}. The second approach classifies ``vanilla" TQFTs, i.e. fully extended and without the extra structures/conditions of the above examples, using the Cobordism Hypothesis \cite{BaezDolan, LurieCob}. This classification is in terms of fully dualizable objects in a higher category. To bridge the two approaches, we must answer the questions: 
\begin{center}\textit{ Can the Cobordism Hypothesis recover the interesting, hand-built examples we know? \\ 
If so, do we know what are the relevant dualizable objects? }\end{center}  
There is evidence that the answer is yes if one allows for relative and non-compact versions of the Cobordism Hypothesis. The relevant dualizable objects were predicted to be those induced by the unit inclusion mentioned in this paper’s title. This follows ideas of Walker, Freed and Teleman in the semisimple case, and was predicted to extend to the non-semisimple case by Jordan and Safronov. The whole story has mainly been communicated in talks; the only written references we are aware of are Walker's notes \cite{WalkerNotes} and Freed's slides \cite{FreedSlides}.

An obstacle to obtaining the Witten--Reshetikhin--Turaev theories from the Cobordism Hypothesis is that these theories are defined on a category of \emph{cobordisms equipped with some extra structure}, which morally comes from the data of a bounding higher manifold. It was noticed by Walker and Freed--Teleman 
that this extra data is actually obtained from the Crane--Yetter 4-TQFT on the bounding manifold. Therefore the WRT theory should be thought of as a boundary theory for the Crane--Yetter theory. We do not know of a formal proof of this statement. An adequate description of relative field theory was given by Freed and Teleman \cite{FreedTeleman} and formalized by Johnson-Freyd and Scheimbauer \cite{JFS}.  

Another obstacle is that the non-semisimple variants are defined on a \emph{restricted class of cobordisms}, namely every 3-cobordism must have non-empty incoming boundary in every connected component. We need to use a non-compact version of the Cobordism Hypothesis to work with this restricted category of cobordisms. This non-compact version appears as an intermediate step in the sketch of proof of the Cobordism Hypothesis proposed by Hopkins and Lurie \cite{LurieCob}. Note that there is independent work in progress of Reutter--Walker and Schommer-Pries in this direction.

A final obstacle is that WRT theories are \emph{not fully extended}. It is known that they extend to the circle, but work of Douglas, Schommer-Pries and Snyder \cite{DSPS}, see also \cite{FreedTelemanGapped}, shows that they extend to the point if and only if they are of Turaev--Viro type. This can be explained by the fact that they come from a relative setting, namely are defined on a category of cobordisms equipped with a bounding higher manifold which we call filled cobordisms, and the point cannot be equipped with a bounding 1-manifold. 

Summing up, one should be able to recover the WRT theories (resp. their non-semisimple variants) from a 4-TQFT and a (resp. non-compact) boundary theory for this 4-TQFT, both of which are fully extended and obtained from the Cobordism Hypothesis. It was proposed by Freed and Teleman in the semisimple case that the 4-TQFT is induced by the modular tensor category $\Vv$, and the boundary theory by the inclusion if the unit in $\Vv$, see the last slide of \cite{FreedSlides}. It was proven in \cite{BJSS} that a possibly non-semisimple modular tensor category $\VV=Ind(\Vv)$ is indeed 4-dualizable in the even higher Morita 4-category of braided tensor categories $\BrTens$, and therefore induces a 4-TQFT under the Cobordism Hypothesis. It was conjectured by Jordan and Safronov in 2019 that this 4-TQFT together with the relative theory induced by the unit inclusion will also recover the non-semisimple variants of WRT
.

This paper gives the first step towards executing the above program. 
We use the framework of \cite{JFS} to prove that the unit inclusion is 3-dualizable (resp. non-compact-3-dualizable), and therefore induces a (resp. non-compact) relative 3-TQFT under the Cobordism Hypothesis. 
In the last section we explain how one can obtain a theory defined on filled cobordisms from a relative theory together with its bulk theory. We state the conjectures that these recover the WRT theories and their non-semisimple variants. Proving these conjectures would answer both questions above in the affirmative.

\subsection*{Context}
The cobordism hypothesis formulated in \cite{BaezDolan}, see \cite{LurieCob} for a sketch of proof, provides a new angle to study and construct Topological Quantum Field Theories. One simply has to find a fully dualizable object in a higher category, and it induces a framed fully extended TQFT. We will study here one particular example of target category, the 4-category $\BrTens$ of braided tensor categories and bimodules between them, or more precisely the even higher Morita category $\Alg_2(\Pr)$ of $\mathbb E_2$-algebras in some 2-category of cocomplete categories.

Even higher Morita categories are defined in \cite{JFS}, see also \cite{Scheimbauer} and \cite{Haugseng}. They form an $(n+k)$-fold Segal space, which we will abbreviate $(n+k)$-category, $\Alg_n(\mathcal{S})$ for $\mathcal{S}$ a $k$-category. It is shown in \cite{GS} that every object in $\Alg_n(\mathcal{S})$ is $n$-dualizable. We study the case $\mathcal{S}=\Pr$, the 2-category locally presentable $\FK$-linear categories, cocontinuous functors and natural transformations, over a field $\FK$ of characteristic zero. 
 It was shown in \cite[Theorem 5.21]{BJS} that fusion categories provide a family of fully dualizable objects in $\BrTens$\footnote{This is where characteristic zero assumption is needed}, and later in \cite[Theorem 1.1]{BJSS} that possibly non-semisimple modular tensor categories are invertible, and hence also fully dualizable. 
 Provided that we can endow these objects with the extra orientation structure needed for the oriented Cobordism Hypothesis, the TQFTs induced by fusion categories are expected to coincide with the Crane--Yetter theories. The TQFTs induced by non-semisimple modular tensor categories are expected to coincide with the ones constructed in \cite{CGHP}.
 
There is another version of the Cobordism Hypothesis for relative field theories. Actually, there are multiple versions as there are multiple notions of relative TQFT. Lurie proposed a definition of a domain wall based on a category of bipartite cobordisms, see \cite[Example 4.3.23]{LurieCob}, and proves the relative Cobordism Hypothesis under the assumption that the ambient category has duals. Stolz and Teichner \cite{StolzTeichner} define a notion of twisted quantum field theories in the context of topological algebras. Freed and Teleman \cite{FreedTeleman} describe relative $n$-TQFTs as morphisms between truncations of fully extended, ideally invertible, $(n+1)$-TQFTs. Johnson-Freyd and Scheimbauer \cite{JFS} give an explicit definition of what a morphism between TQFTs is and exhibit three different notions of strong-, lax- and oplax-twisted quantum field theories, which also makes sense when the bulk theories are only $(n+\varepsilon)$-TQFTs. We will be mostly interested in their notion of oplax-twisted quantum field theory.

We consider the category $\BrTens^\to$ of arrows in our chosen target category, where morphisms are ``oplax" squares filled by a 2-morphism. There is a well-defined notion of source and target for objects and morphisms in this arrow category. Given $\ZZ$ a fully extended $(n+\varepsilon)$-TQFT, a relative theory to $\ZZ$ is a symmetric monoidal functor $\RR:\Bord_n\to \BrTens^\to$ whose source is trivial and whose target coincides with $\ZZ$. The Cobordism Hypothesis applies directly in this context, namely $\RR$ can be reconstructed from its value on the point $\RR(pt): \unit \to \ZZ(pt)$ which has to be fully dualizable in $\BrTens^\to$. We say that the 1-morphism $\RR(pt)$ has to be fully oplax-dualizable.

In this paper we study the dualizability of the 1-morphism induced by the inclusion of the unit $\eta: \Vect_\FK \to \VV$. The braided monoidal functor $\eta$ induces a $\Vect_\FK$-$\VV$-central algebra $\AA_\eta$ which is $\VV$ as a tensor category with bimodule structure induced by $\eta$. When it is oplax-dualizable, it induces a relative TQFT to the one of $\VV$. 
There is a stronger notion of dualizability for a 1-morphism needed to induced a domain wall in the sense of Lurie. 
It is already known that the 1-morphism $\AA_\eta$ is always 1-dualizable in the strong sense, by \cite{GS}, 2-dualizable as soon as $\VV$ is cp-rigid, and 3-dualizable as soon as $\VV$ is fusion, by \cite{BJS}. 

We show that fusion is a criterion for 3-dualizability, but not for 3-oplax-dualizability, emphasizing the difference between these notions. We study oplax-dualizability in detail, including the non-semisimple cases. It is expected that the induced oplax-twisted theory corresponds to the Witten--Reshetikhin--Turaev TQFT seen as relative to the Crane--Yetter 4-TQFT in Walker's \cite{WalkerNotes} and Freed--Teleman's \cite{FreedSlides} picture.

In the modular non-semisimple case, we will show that $\AA_\eta$ is \emph{not} 3-oplax-dualizable. We can only hope for a partial dualizability, and a partially defined TQFT. It turns out that Lurie's sketch of proof of the cobordism hypothesis is building a TQFT inductively from the dualizability data of the value on the point, and if this object is not fully dualizable, this process will simply stop partway through. Note again related ongoing work of Reutter--Walker and Schommer-Pries. In our case this process will stop before the very last step, and we obtain a theory defined on cobordisms with outgoing boundary in every connected component, which we call a non-compact TQFT, see Section \ref{subsectNonCompact}. 

\subsection*{Results}
The unit inclusion in a braided tensor category $\VV\in \BrTens$ gives a braided monoidal functor $\eta: \Vect_\FK\to \VV$. We work over a field $\FK$ of characteristic zero. Using Definition \ref{defAF}, it induces a $\Vect_\FK$-$\VV$-central algebra $\AA_\eta$, i.e. a 1-morphism in $\BrTens$. When we see this morphism as an object of the oplax arrow category $\BrTens^\to$, we denote it $\Aflat_\eta$. We are interested in the adjunctibility of $\AA_\eta$ and in its oplax-dualizability, i.e. in the dualizability of $\Aflat_\eta$.

We recall Lurie's sketch of proof of the non-compact cobordism hypothesis and introduce the corresponding notion of non-compact-$n$-dualizable object in Section \ref{subsectNonCompact}. We characterize oplax-dualizability of the unit inclusion:
\begin{myth}[Theorems \ref{thmVcpgen}, \ref{thmCritnc3d} and \ref{thmCrit3oplax}]\label{thmMain} Let $\VV\in \BrTens$ be a braided tensor category, and $\Aflat_\eta$ the object of $\BrTens^\to$ induced by the inclusion of the unit. Then:
\begin{enumerate}
    \item $\Aflat_\eta$ is always 2-dualizable.
\end{enumerate}
If $\VV$ has enough compact-projectives, then:
\begin{enumerate} \setcounter{enumi}{1}
    \item $\Aflat_\eta$ is non-compact-3-dualizable if and only if $\VV$ is cp-rigid.
    \item $\Aflat_\eta$ is 3-dualizable if and only if $\VV$ is the free cocompletion of a small rigid braided monoidal category if and only if $\VV$ is cp-rigid with compact-projective unit.
\end{enumerate}
\end{myth} 
In particular, on the examples of interest for Section \ref{sectAppli}:
\begin{corr} \label{corrModular}
Let $\Vv$ be a modular tensor category in the sense of \cite{DGGPR}, $\VV := Ind(\Vv)$ its Ind-completion, and $\Aflat_\eta$ induced by the unit inclusion in $\VV$. Then: \begin{enumerate}
    \item If $\Vv$ is semisimple, $\Aflat_\eta$ is 3-dualizable and induces a 3-TQFT $\RR_\VV$ with values in $\BrTens^\to$.
    \item If $\Vv$ is non-semisimple, $\Aflat_\eta$ is not 3-dualizable, but is non-compact-3-dualizable and induces a non-compact-3-TQFT $\RR_\VV$ with values in $\BrTens^\to$.
\end{enumerate}
\end{corr}
We describe the dualizability data of $\Aflat_\eta$ explicitly, which gives the values of $\RR_\VV$ on elementary handles. In dimension 2, the handle of index 2 is mapped to some mate of the unit inclusion $\eta$. The handle of index 1 is mapped to some mate of the tensor product $T$. And the handle of index 0 is mapped to some mate of the ``balanced tensor product" $T_{bal}:\VV \underset{\VV\boxtimes\VV}{\boxtimes}\VV \to \VV$ which is induced by $T$ on the relative tensor product.

To determine this dualizability data, we use the fact that the dualizability of $\Aflat_\eta$ is equivalent to the dualizability of $\VV$ and the right-adjunctibility of $\AA_\eta$, see \cite[Theorem 7.6]{JFS} and Section \ref{subsectOplax}. We remark that the adjunctibility of $\AA_\eta$ often implies the relevant dualizability of $\VV$, see Remark \ref{rmkAdjImpliesDual}, which is a priori a phenomenon specific to the unit inclusion. 

Freed and Teleman study the dualizability of the unit inclusion in the 3-category $\Alg_1(Rex_\Cc)$ in \cite[Theorem B]{FreedTelemanGapped}. Here $Rex_\Cc$ is the 2-category of finitely cocomplete categories and finitely cocontinuous functors. They show that $\Vv \in \Alg_1(Rex_\Cc)$ is finite rigid semisimple if and only if $\MM_\eta$ is 2-dualizable, i.e. lies in a subcategory with duals. The forward implication is \cite{DSPS}. We can give an analogous statement one categorical number higher:
\begin{myth}[Theorem \ref{thmFTplus1}]
Suppose $\VV\in \BrTens$ has enough compact-projectives. Then $\AA_\eta$ is 3-dualizable if and only if $\VV$ is finite rigid semisimple.
\end{myth}
Note that the full dualizability of $\AA_\eta$ is indeed stronger that its oplax-dualizability. They are however expected to be equivalent when $\VV$ is itself fully dualizable, see Remark \ref{rmkWillStewart}. This actually implies some non-dualizability results on $\VV$, see Remark \ref{rmkNon4DV}

We study in Section \ref{sectFunBimod} the dualizability of the 1- (resp. 2-) morphism induced by a braided monoidal (resp. bimodule monoidal) functor $F$, which we denote $\AA_F$ (resp. $\MM_F$). This allows us to give explicitly the adjunctibility data of $\AA_\eta$, see Figure \ref{fig:DDtreeCprigid}. We will use the notation
\begin{center}
    \begin{tikzpicture}[baseline = 10pt, xscale = 1.2]
\node (A) at (1,1) {morphism};
\node (L) at (0,0) {counit $\vert$ unit};
\node (R) at (2,0) {counit $\vert$ unit};
\draw[->] (A) -- (R) node[near end,above right=-3pt, scale = 0.8]{Right adjoint};
\draw[->] (A) -- (L) node[near end,above left=-3pt, scale = 0.8]{Left adjoint};
\end{tikzpicture}
\end{center} 
to depict adjunctibility data. We will give a further example explaining this notation in Figure \ref{fig:cupDD}. 
\begin{figure}[h]
\centering
\begin{tikzpicture}[baseline = 10pt, xscale = 1.4, yscale = 1.2]
\node (A) at (2,2) {$\AA_\eta$};
\node[inner sep = 20pt] (L) at (0,1) {$\vert$};
\node[inner sep = 20pt] (R) at (4,1) {$\vert$};
\node[left of= R] (Rcounit)  {$\MM_T$};
\node[right of= R] (Runit){$\MM_\eta$};
\node[left of= L] (Lcounit) {$\overline{\MM}_\eta$};
\node[right of= L] (Lunit) {$\overline{\MM}_T$};
\node (LLc) at (-2,0) {$T\vert \eta$};
\node[gray] (RLc) at (-0.5,0) {$\widetilde{\eta^R}\vert T^R$};
\node (LLu) at (0.5,0) {$T_{bal}\vert T$};
\node[gray] (RLu) at (1.5,0) {$T^R\vert {T_{bal}^R}$};
\node[gray] (LRc) at (2.5,0) {$T^R\vert {T_{bal}^R}$};
\node (RRc) at (3.5,0) {$T_{bal}\vert T$};
\node[gray] (LRu) at (4.5,0) {$\widetilde{\eta^R}\vert T^R$};
\node (RRu) at (6,0) {$T\vert \eta$};
\draw[->] (A) -- (R) node[near end,above right]{\footnotesize $\overline{\AA}_\eta$};
\draw[->] (A) -- (L) node[near end,above left]{\footnotesize $\overline{\AA}_\eta$};
\draw[->, gray] (Rcounit) -- (LRc) node[near end,above left, gray]{\tiny $\overline{\MM}_T$};
\draw[->] (Rcounit) -- (RRc) node[near end,above right]{\tiny $\overline{\MM}_T$};
\draw[->, gray] (Runit) -- (LRu) node[near end,above left, gray]{\tiny $\overline{\MM}_\eta$};
\draw[->] (Runit) -- (RRu) node[near end,above right]{\tiny $\overline{\MM}_\eta$};
\draw[->] (Lcounit) -- (LLc) node[near end,above left]{\tiny $\MM_\eta$};
\draw[->, gray] (Lcounit) -- (RLc) node[near end,above right, gray]{\tiny $\MM_\eta$};
\draw[->] (Lunit) -- (LLu) node[near end,above left]{\tiny $\MM_T$};
\draw[->, gray] (Lunit) -- (RLu) node[near end,above right, gray]{\tiny $\MM_T$};
\end{tikzpicture}
\label{fig:DDtreeCprigid}
\caption{Adjunctibility data of the unit inclusion. The whole description (including gray) holds for $\VV$ cp-rigid, see Proposition \ref{propVcprigid}, and the black subset holds when $\VV$ has enough compact-projectives, see Theorem \ref{thmVcpgen}. The functor $\widetilde{\eta^R}$ is the essentially unique cocontinuous functor that agrees with $\eta^R$ on the compact-projectives.}
\end{figure}

We studied the unit inclusion, but a version of our arguments still work for any bimodule induced by a braided monoidal functor. Instead of a necessary and sufficient condition, we only have a sufficient condition:
\begin{myth}[Theorem \ref{thm3oplaxAnyBrMonFun}]
Let $F: \VV \to \WW$ be a braided monoidal functor between two objects of $\BrTens$. Then the object $\Aflat_F \in \BrTens^\to$ induced by $F$ is 2-dualizable. It is non-compact-3-dualizable as soon as $\VV$ and $\WW$ are cp-rigid. In this case, it is 3-dualizable if and only if $F$ preserves compact-projectives.
\end{myth}

\subsection*{Applications}
Non-semisimple variants of Witten--Reshetikhin--Turaev TQFTs were introduced in \cite{BCGP} and \cite{DGGPR}. They are defined on a restricted class of decorated cobordisms, including in particular cobordisms with incoming boundary in every connected component. This matches the notion of non-compact TQFT from Lurie, up to orientation reversal, and this is the part of the theory that we expect to obtain. We will focus on the TQFTs from \cite{DGGPR} here, and actually their extension to the circle by \cite{DeRenziDGGPR}. They are defined for possibly non-semisimple modular categories, whose Ind-completions have been found to be 4-dualizable, and actually invertible, in $\BrTens$ by \cite{BJSS}. Using the Cobordism Hypothesis, for $\Vv$ a modular tensor category and $\VV:=Ind(\Vv)$ its Ind-completion, there is an essentially unique framed 4-TQFT $\ZZ_\VV:\Bord_4^{fr}\to \BrTens$ with $\ZZ_\VV(pt)=\VV$.

We know from Corollary \ref{corrModular} that if $\VV$ is semisimple (resp. non-semisimple) the unit inclusion is 3-oplax-dualizable (resp. non-compact-3-oplax-dualizable). Using the Cobordism Hypothesis, it induces a framed (resp. non-compact) 3-TQFT $\RR_\VV:\Bord_3^{fr}\to \BrTens^\to$ (resp. $\RR_\VV:\Bord_3^{fr, nc}\to \BrTens^\to$) relative to $\ZZ_\VV$.
We give conjectures that these theories can be oriented.
\begin{conjecture}[Conjectures \ref{conjSO3} and \ref{conjSO45}] Let $\Vv$ be a modular tensor category and $\VV= \Ind(\Vv)$, which is invertible and in particular 5-dualizable by \cite{BJSS}. Then,
\begin{enumerate}
    \item The ribbon structure of $\mathscr{V}$ induces an $SO(3)$-homotopy-fixed-point structure on $\VV$.
    \item The ribbon structure of $\eta$ induces an $SO(3)$-homotopy-fixed-point structure on $\Aflat_\eta$.
    \item A choice of modified trace $\mt$ on $\mathscr{V}$ induces an $SO(4)$-homotopy-fixed-point structure \mbox{on $\VV$}.
\end{enumerate}
Given a modified trace $\mt$, let $d(\Vv)_\mt$ denote the global dimension of $\Vv$ computed using $\mt$, defined as the value on $S^4$ of the $(3+1)$-TQFT of \cite{CGHP} with the same input.
\begin{enumerate} 
\setcounter{enumi}{3}
    \item Exactly two modified traces induce $SO(5)$-homotopy-fixed-point structures on $\VV$, namely $\pm\mathscr D_\mt^{-1}\mt$ for $\mathscr D_\mt$ a square root of the global dimension $d(\Vv)_\mt$.
\end{enumerate}
\end{conjecture}
Let us include here the conjecture that the TQFTs of \cite{CGHP} compute the (3+1)-part of the fully extended TQFT associated with $\Vv$, which we can state now that we have conjectured orientation structures.
\begin{conjecture}\label{conjCGHP} Let $\Vv$ be a modular tensor category. Choose $\mt$ a modified trace on $\Vv$ and let $\ZZ_\VV$ be the associated oriented 4-TQFT. Then one has a natural isomorphism $$\mathscr S_\Vv \simeq h_1\Omega^{3}\ZZ_\VV$$ between the (3+1)-TQFT defined in \cite{CGHP} and the (3+1)-part of $\ZZ_\VV$.
\end{conjecture}

We construct the ``anomalous" theory $\AA_\VV:\Bord_3^{filled}\to \Tens$ associated with $\RR_\VV$ and $\ZZ_\VV$. It is defined on a 3-category of cobordisms equipped with a filling, i.e. a bounding higher manifold, which degenerates to the more usual $\widetilde \cob$ on which WRT-type theories are defined. The anomalous theory is non-compact when $\RR_\VV$ is.

We can now state the main conjectures. We claim that one can recover WRT and DGGPR theories from the cobordism hypothesis using the construction we described: 
\begin{conjecture}[Conjecture \ref{conjSS}]
Let $\mathscr{V}$ be a semisimple modular tensor category with a chosen square root of its global dimension. The anomalous theory $\AA_\VV$ induced by the associated oriented 4-TQFT $\ZZ_\VV$ and oriented oplax-$\ZZ_\VV$-twisted 3-TQFT $\RR_\VV$ recovers the once-extended Witten--Reshetikhin--Turaev theory as its 321-part.
\end{conjecture}
\begin{conjecture}[Conjecture \ref{conjNonSS}]
Let $\mathscr{V}$ be a non-semisimple modular tensor category with a chosen modified trace $\mt$ and square root of its global dimension. The anomalous theory $\AA_\VV$ induced by the associated oriented 4-TQFT $\ZZ_\VV$ and oriented non-compact oplax-$\ZZ_\VV$-twisted 3-TQFT $\RR_\VV$ recovers the once-extended DGGPR theory for cobordisms with trivial decoration as its 321-part.
\end{conjecture} 
We show that the values on the circle coincide by computing $\RR_\VV$ explicitly and using the factorization homology description of $\ZZ_\VV$. For higher dimensions, one needs to identify the values of $\ZZ_\VV$ on 3-manifolds with skein modules, which is at this day still conjectural \cite[Conjecture 9.10]{JohnsonFreydHeisenbergPicture}.

An interesting consequence of these conjectures is that since the WRT and DGGPR theories only correspond to the 321-part of $\AA_\VV$, they actually extend down. They do not descend to the point, which is not null-bordant, but to the pair of points $S^0$.

Note that we use both the oplax-twisted 3-TQFT and the 4-TQFT in this construction. Therefore, not every case of 3-oplax-dualizability in Theorem \ref{thmMain} induces such a theory. We also need $\VV$ to be 4-dualizable. The assumption that $\ZZ_\VV$ is invertible however can be dropped. The anomalous theory would then strongly depend on the filling, and give interesting invariants of 4-manifolds with boundary.

The construction of the anomalous theory $\AA_\VV$ using a bounding manifold is needed to recover the usual constructions of WRT and DGGPR theories. It is also necessary for some applications, e.g. to obtain a scalar invariant of 3-manifold. However, one could argue that the more fundamental object is the fully extended twisted 3-TQFT $\RR_\VV$. It does not assign a scalar to a 3-manifold, but an element in a one-dimensional vector space: the state space of the invertible Crane--Yetter TQFT. If one is happy to allow this feature, then one can argue that WRT is a fully extended theory, with values in the oplax arrow category $\BrTens^\to$.

\subsection*{Acknowledgements}
I would like to thank my advisors Francesco Costantino and David Jordan for their invaluable help and encouragements. I am grateful to Pavel Safronov, Will Stewart and Marco De Renzi for many helpful discussions, to my PhD referees Anna Beliakova and Claudia Scheimbauer and to the anonymous referee for many useful remarks on a first version of this article. This research is part of my PhD thesis and took place both at the Universit\'e Toulouse 3 Paul Sabatier and at the University of Edinburgh.

\section{Relative and Non-compact TQFTs} \label{sectDual} In this paper we will study the dualizability of a 1-morphism. What exact kind of dualizability we are interested in is dictated by the relative cobordism hypothesis: we want a 1-morphism that will induce a relative TQFT. It turns out that there are multiple notions of relative TQFTs, and therefore multiple interesting notions of dualizability for a 1-morphism.

Throughout, we will use the expression $n$-category to mean $(\infty,n)$-category, and more precisely complete $n$-fold Segal space. For $j\geq k$, we write $\circ_k$ for the composition of $j$-morphisms in the direction of $k$-morphisms. We write $\Id^k_f$ for taking $k$-times the identity of $f$.

\subsection{Review of relative TQFTs}\label{subsectRelTQFT} We recall the notions of relative TQFTs that will be our motivation. Let $\CC$ be a symmetric monoidal $n$-category. We distinguish two flavors.

The first is purely topological. Lurie defines a category $\Bord_n^{dw}$ of bipartite cobordisms with two different colors for the bulk and interfaces between them, see \cite[Example 4.3.23]{LurieCob}. There are in particular manifolds with only one color and without interfaces. This induces two inclusions $\Bord_n \to \Bord_n^{dw}$.
\begin{mydef}[Lurie]
A domain wall between two theories $\ZZ_1, \ZZ_2: \Bord_n \to \CC$ is a symmetric monoidal functor $\Bord_n^{dw}\to \CC$ that restricts to $\ZZ_1$ and $\ZZ_2$ on manifolds with one color.
\end{mydef}
In particular, the interval with an interface point in the middle induces a morphism $\ZZ_1(pt) \to \ZZ_2(pt)$. Freed and Teleman describe a notion of relative TQFT by means of such morphisms for every values of $\ZZ_1$ and $\ZZ_2$ on manifolds of dimension strictly less than $n$, see \cite{FreedTeleman}. They mention that their notion should be equivalent to Lurie's notion of domain wall. A more detailed comparison will appear in William Stewart's PhD thesis.

The second notion focuses on the algebraic flavour of Freed--Teleman's description. One can drop the assumption that $\ZZ_1$ and $\ZZ_2$ are well defined on $n$-manifolds because these don't appear. Johnson-Freyd and Scheimbauer define three different notions of an $n$-category of arrows in an $n$-category. We will focus on the oplax one $\CC^\to$.
\begin{mydef}[sketch, see Definition 5.14 in \cite{JFS}]
Let $\CC$ be a symmetric monoidal $n$-category. The \emph{symmetric monoidal $n$-category $\CC^\to$ of oplax arrows in $\CC$} has:\vskip8pt \noindent
\begin{tabular}{p{\widthof{$k$-morphisms}}lp{\textwidth-\widthof{1-morphisms : \quad\quad\ }}}
objects &:& triples $f=(s_f,t_f,f^\#)$ where $s_f$ and $t_f$ are objects of $\CC$ and $f^\#: s_f \to t_f$ is a 1-morphism\\
1-morphisms \newline $f \to g$ &:& triples $h=(s_h,t_h,h^\#)$ where $s_h: s_f \to s_g$ and $t_h:t_f\to t_g$ are 1-morphisms, and $h^\#: g^\# \circ s_h \Rightarrow t_h \circ f^\#$ is a 2-morphism\\
& & $\quad\quad\vdots$\\
$k$-morphisms\newline $a\to b$  &:& triples $f= (s_f, t_f, f^\#)$ where $s_f: s_a \to s_b$ and $t_f: t_a \to t_b$ are $k$-morphisms in $\CC$, and $f^\#$ is a $k+1$-morphism in $\CC$ from the composition of some whiskerings of $b^\#$ and $s_f$ to the composition of some whiskerings of $t_f$ and $a^\#$.
\end{tabular}\vskip8pt
It has two symmetric monoidal functors $s,t:\CC^\to \to \CC$.

To avoid confusion, when we see a 1-morphism $f$ of $\CC$ as an object of $\CC^\to$ we will denote it $f^\flat$. The notation comes from $(f^\flat)^\#=f$.
\end{mydef}
\begin{mydef}[Definition 5.16 in \cite{JFS}]
Let $\CC$ be a symmetric monoidal $n$-category and $\ZZ_1, \ZZ_2: \Bord_{n-1} \to \CC$ two categorified $(n-1)$-TQFTs. An \emph{oplax-$\ZZ_1$-$\ZZ_2$-twisted $(n-1)$-TQFT} is a symmetric monoidal functor $$\RR: \Bord_{n-1}\to \CC^\to$$ such that $s(\RR)= \ZZ_1$ and $t(\RR)= \ZZ_2$.
\end{mydef}
The name and strategy come from \cite{StolzTeichner}.

We will use the formalism of Johnson-Freyd and Scheimbauer in this paper. For application, see Section \ref{sectAppli}, we are interested in the case where $\ZZ_1$ is the trivial theory and $\ZZ_2$ is well defined on $n$-manifolds. If $\ZZ: \Bord_n \to \CC$ is defined on $n$-manifolds, we will say oplax-$\ZZ$-twisted theory for oplax-$\Triv$-$\ZZ\vert_{\Bord_{n-1}}$-twisted theory. 
Under this extra hypothesis, which was made in \cite{FreedTeleman}, the notion of oplax-twisted field theory should agree with Lurie's notion of domain walls, see Remark \ref{rmkWillStewart}.

\subsection{Dualizability data}\label{subsectDual}
Let us first recall the multiple notions of dualizability and adjunctibility for morphisms in a symmetric monoidal $n$-category $\CC$.
\subsubsection{Definitions}
\begin{mydef}
Let $\CC$ be a bicategory, and $f : x \to y$ a 1-morphism in $\CC$. A right adjoint for $f$ is a morphism $f^R:y\to x$ together with two 2-morphisms $\varepsilon: f\circ f^R \Rightarrow \Id_y$ called the counit and $\eta: \Id_x \Rightarrow f^R \circ f$ called the unit, satisfying the so-called snake identities:\\
\begin{tikzpicture}[xscale = 2, baseline = 1.3cm]
\node (sIdx) at (0,3) {$x$};
\node (tIdx) at (0,1.5) {$x$};
\draw[double] (sIdx)--(tIdx) node[midway, left]{$\Id_x$};
\node (tf1) at (0,0) {$y$};
\draw[->] (tIdx) -- (tf1) node[midway, left]{$f$};
\node (seta) at (1,3) {$x$};
\node (meta) at (1,2) {$y$};
\node (teta) at (1,1) {$x$};
\draw[->] (seta) -- (meta) node[midway, left]{$f$};
\draw[->] (meta) -- (teta) node[midway, left]{$f^R$};
\node (tf2) at (1,0) {$y$};
\draw [->] (teta) -- (tf2) node[midway, left]{$f$};
\node (sf) at (2,3){$x$};
\node (sIdy) at (2,1.5){$y$};
\draw[->] (sf) -- (sIdy) node[midway, right]{$f$};
\node (tIdy) at (2,0){$y$};
\draw[double] (sIdy)--(tIdy) node[midway, right]{$\Id_y$};
\draw[very thin, gray!50] (sIdx) -- (seta);
\draw[very thin, gray!50] (tIdx) -- (teta);
\draw[very thin, gray!50] (tf1) -- (tf2);
\draw[very thin, gray!50] (seta) -- (sf);
\draw[very thin, gray!50] (meta) -- (sIdy);
\draw[very thin, gray!50] (tf2) -- (tIdy);
\draw[double, ->] (0.2, 2.2) -- (0.75,2) node[midway, above]{$\eta$};
\draw[double, ->] (1.2, 0.8) -- (1.8,0.6) node[midway, above]{$\varepsilon$};
\node at (1.5,2.5){$=$}; \node at (0.5,0.5){$=$};
\node (Fsf) at (-1,2.3) {$x$};
\node (Ftf) at (-1,0.7) {$y$};
\draw [->] (Fsf) -- (Ftf) node[midway, left]{$f$};
\node (Lsf) at (3,2.3) {$x$};
\node (Ltf) at (3,0.7) {$y$};
\draw [->] (Lsf) -- (Ltf) node[midway, right]{$f$};
\draw[very thin, gray!50] (Fsf) -- (sIdx);
\draw[very thin, gray!50] (Ftf) -- (tf1);
\draw[very thin, gray!50] (sf) -- (Lsf);
\draw[very thin, gray!50] (tIdy) -- (Ltf);
\node at (-0.5,1.5){$\simeq$};
\node at (2.5,1.5){$\simeq$};
\end{tikzpicture} $\quad=\quad \Id_f$ \quad 
and,\\
\begin{tikzpicture}[xscale = 2, baseline = 1.3cm]
\node (sfR1) at (0,3) {$y$};
\node (sIdx) at (0,1.5) {$x$};
\node (tIdx) at (0,0) {$x$};
\draw[double] (sIdx)--(tIdx) node[midway, left]{$\Id_y$};
\draw[->] (sfR1) -- (sIdx) node[midway, left]{$f^R$};
\node (seta) at (1,2) {$x$};
\node (meta) at (1,1) {$y$};
\node (teta) at (1,0) {$x$};
\node (sfR2) at (1,3) {$y$};
\draw[->] (seta) -- (meta) node[midway, left]{$f$};
\draw[->] (meta) -- (teta) node[midway, left]{$f^R$};
\draw [->] (sfR2) -- (seta) node[midway, left]{$f^R$};
\node (sIdy) at (2,3){$y$};
\node (tIdy) at (2,1.5){$y$};
\node (tfR) at (2,0){$x$};
\draw[->] (tIdy) -- (tfR) node[midway, right]{$f^R$};
\draw[double] (sIdy)--(tIdy) node[midway, right]{$\Id_y$};
\draw[very thin, gray!50] (sIdx) -- (seta);
\draw[very thin, gray!50] (tIdx) -- (teta);
\draw[very thin, gray!50] (sfR1) -- (sfR2);
\draw[very thin, gray!50] (sfR2) -- (sIdy);
\draw[very thin, gray!50] (meta) -- (tIdy);
\draw[very thin, gray!50] (teta) -- (tfR);
\draw[double, ->] (0.2, 0.7) -- (0.75,1) node[midway, above]{$\eta$};
\draw[double, ->] (1.2, 1.9) -- (1.8,2.2) node[midway, above]{$\varepsilon$};
\node at (1.5,0.5){$=$}; \node at (0.5,2.5){$=$};
\node (Fsf) at (-1,2.3) {$y$};
\node (Ftf) at (-1,0.7) {$x$};
\draw [->] (Fsf) -- (Ftf) node[midway, left]{$f^R$};
\node (Lsf) at (3,2.3) {$y$};
\node (Ltf) at (3,0.7) {$x$};
\draw [->] (Lsf) -- (Ltf) node[midway, right]{$f^R$};
\draw[very thin, gray!50] (Fsf) -- (sfR1);
\draw[very thin, gray!50] (Ftf) -- (tf1);
\draw[very thin, gray!50] (sIdy) -- (Lsf);
\draw[very thin, gray!50] (tfR) -- (Ltf);
\node at (-0.5,1.5){$\simeq$};
\node at (2.5,1.5){$\simeq$};
\end{tikzpicture} $\quad=\quad \Id_{f^R}$.\\
We say that $f$ has a right adjoint $f^R$ and that $f^R$ has a left adjoint $f$.

This definition extends to higher categories. Let $\CC$ be an $n$-category, $2\leq k\leq n$ and $f:x\to y$ a $k$-morphism between two $k-1$-morphisms $x,y: a\to b$ in $\CC$. A right adjoint for $f$ is a right adjoint for $f$ seen as a 1-morphism in the bicategory $h_2(\Hom_\CC(a,b))$. If $k=1$ then we demand a right adjoint of $f$ in $h_2(\CC)$. If $k=0$ and $\CC$ is a monoidal category then we demand a right adjoint of $f$ in $h_2(B\CC)$, where $B\CC$ is the one object $n+1$-category with endomorphisms of the object being $\CC$, and composition the monoidal structure of $\CC$, i.e. $X \circ_0 Y := X \otimes Y$.
\end{mydef}
\begin{mydef}
Let $\CC$ be a symmetric monoidal $n$-category. It is said to have duals up to level $m$ if every $k$-morphism of $\CC$, $0\leq k<m$, has both a left and a right adjoint. It is said to have duals if it has duals up to level $n$.

An object $X \in \CC$ is said $m$-dualizable if it lies in a sub-$n$-category with duals up to level $m$. It is called fully dualizable if it is $n$-dualizable.
\end{mydef}

There are multiple notions of dualizability, or adjunctibility, for morphisms and higher morphisms in $\CC$. Following \cite{LurieCob} one defines:
\begin{mydef}
 A $k$-morphism $f$ of $\CC$ is said\emph{ $m$-dualizable} if it lies in a sub-$n$-category with duals up to level $m+k$. It is called \emph{fully dualizable} if it is $(n-k)$-dualizable.
\end{mydef}
Following \cite{JFS} one gets a few more notions. For simplicity we focus on 1-morphisms.
\begin{mydef}
A $1$-morphism $f:X \to Y$ of $\CC$ is said \emph{$m$-oplax-dualizable} if it is $m$-dualizable as an object $f^\flat$ of $\CC^\to$. It is said \emph{$m$-lax-dualizable} if it is $m$-dualizable as an object of $\CC^\downarrow$, where $\CC^\downarrow$ is the category of lax arrows defined in \cite[Definition 5.14]{JFS}. 
\end{mydef}
\begin{mydef}
A $k$-morphism $f$ is said to be \emph{left} (resp. \emph{right}) \emph{adjunctible} if it has a left (resp. right) adjoint, and \emph{adjunctible} if it has arbitrary left and right adjoints ($(f^L)^L$, $(f^R)^R$ and so on...). It is said to be \emph{$m$-times} (resp. \emph{left}, \emph{right}) \emph{adjunctible} if it is $m-1$-times (resp. left, right) adjunctible and every unit/counit witnessing this are themselves (resp. left, right) adjunctible. We sometimes abbreviate $m$-times adjunctible as $m$-adjunctible.
\end{mydef}

Note that being (left, right) adjunctible is only a condition on the morphism while being (lax, oplax) dualizable is also a condition on its source and target.
\begin{myth}[Theorem 7.6 in \cite{JFS}]\label{thmCaracOplaxDualiz}
A 1-morphism $f:X \to Y$ of $\CC$ is $m$-oplax-dualizable if and only if $X$ and $Y$ are both $m$-dualizable and $f$ is $m$-times right adjunctible.

 Similarly, it is $m$-lax-dualizable if and only if $X$ and $Y$ are both $m$-dualizable and $f$ is $m$-times left adjunctible.
\end{myth}
Similarly, a 1-morphism $f:X \to Y$ is $m$-dualizable if and only if it is $m$-times adjunctible and its source and targets are $m+1$-dualizable. 

\subsubsection{Redundancy in the dualizability data} \label{subsectRedundancyDD} The dualizability data of a morphism grows very fast: there are four units/counits for the left and right adjunctions, and this does not consider taking the right adjoint of the right adjoint and so on. In particular, checking $n$-adjunctibility of a morphism seems tedious. It turns out that there is a lot of redundancy in this data, especially if we are only interested in dualizability properties.

Let us begin with some notations. Let $f$ be a $k$-morphism in an $n$-category. We say that $Radj(f)$ (resp. $Ladj(f)$) exists if $f$ has a right (resp. left) adjoint, in which case we denote this adjoint $Radj(f)$ (resp. $Ladj(f)$), and the unit and counit of the adjunction $Ru(f)$ and $Rco(f)$ (resp. $Lu(f)$ and $Lco(f)$). When these adjoints exist we will display the right and left dualizability data as:
\begin{center}
\begin{tikzpicture}[baseline = 10pt, xscale = 1.2]
\node (A) at (1,1) {$f$};
\node (R) at (2,0) {\small $Rco(f)\vert Ru(f)$};
\draw[->] (A) -- (R) node[scale = 0.8, near end,above right=-3pt]{$Radj(f)$};
\end{tikzpicture} \quad\quad and \quad\quad \begin{tikzpicture}[baseline = 10pt, xscale = 1.2]
\node (A) at (1,1) {$f$};
\node (L) at (0,0) {\small $Lco(f)\vert Lu(f)$};
\draw[->] (A) -- (L) node[scale = 0.8, near end,above left=-3pt]{$Ladj(f)$};
\end{tikzpicture}
\end{center}
See Figure \ref{fig:cupDD} for an example. 
Note that adjoints, units and counits are only defined up to isomorphism, and we may write $Radj(f)$, $Ru(f)$, $Rco(f)$ for any choice of adjoint, unit, counit. The fact that these choices do not matter will be shown in point 1 of Proposition \ref{prop:coherenceDualiz}.
\begin{figure}[h]
    \centering
    \begin{tikzpicture}[baseline = 10pt, xscale = 3, yscale = 1.7]
\node (A) at (1,1) {\rcup};
\node (L) at (0,0) {\begin{tikzpicture}[scale = 0.5]
    \draw (0,0) arc(-90:90:1);
    \draw (0,0) arc(-90:-270:1);
    \draw (0,0) arc(-90:90:2 and 1);
\end{tikzpicture} $\quad\vert\quad$
\begin{tikzpicture}[scale = 0.7]
    \draw (0,0)--++(0,1);
    \draw (1,0)--++(0,1);
    \draw[dashed] (1.5,0.5) arc(180:90:0.5 and 0.35);
    \draw (2,0.85) arc(90:0:0.5 and 0.35);
    \draw[dashed] (1.5,1.5) arc (180:270:0.5 and 0.35);
    \draw (2,1.15) arc (270:360:0.5 and 0.35);
    \draw (0,0)--++(0.9,0.3);
    \draw[dashed] (1,0.33)--(1.5,0.5);
    \draw (1,0)--++(1.5,0.5);
    \draw (0,1)--++(1.5,0.5);
    \draw (1,1)--++(1.5,0.5);
    \draw (2,0.85) .. controls (1.7, 0.78).. (1.4, 0.8);
    \draw (2,1.15) .. controls (1.7, 1.05).. (1.25, 0.7);
\end{tikzpicture}};
\node (R) at (2,0) {\begin{tikzpicture}[scale = -0.7, baseline = -15pt]
    \draw (0,0)--++(0,1);
    \draw[dashed] (1,0)--++(0,1);
    \draw (1.5,0.5) arc(180:0:0.5 and 0.35);
    \draw (1.5,1.5) arc (180:360:0.5 and 0.35);
    \draw (0,0)--++(1.5,0.5);
    \draw (1,0)--++(1.5,0.5);
    \draw (0,1)--++(1.5,0.5);
    \draw[dashed] (1,1)--++(0.75,0.25);
    \draw (1.81,1.27)--(2.5,1.5);
    \draw (2,0.85) .. controls (1.7, 0.78).. (1.2, 0.8);
    \draw (2,1.15) .. controls (1.7, 1.05).. (1.05, 0.7);
\end{tikzpicture} $\quad\vert\quad$ 
\begin{tikzpicture}[scale = 0.5]
    \draw (0,0) arc(-90:90:1);
    \draw[dashed] (0,0) arc(-90:-270:1);
    \draw (0,0) arc(-90:-270:2 and 1);
\end{tikzpicture}};
\draw[->] (A) -- (R) node[near end,above right=3pt, scale = 0.8]{\lcap};
\draw[->] (A) -- (L) node[near end,above left=3pt, scale = 0.8]{\lcap};
\end{tikzpicture}
\caption{Example of the notation for the adjunctibility data of the cup $\cup\hskip-0.8pt\hat{}:\emptyset\to  \bullet^-\sqcup\bullet^+$ in the bicategory $\cob_{0,1,2}$ of 0, 1 and 2-dimensional oriented cobordisms}
    \label{fig:cupDD}
\end{figure} 
\begin{mydef}
We say that two $k$-morphisms $f$ and $g$ have same dualizability properties, which we denote $f \simdual g$, if for every finite sequence \begin{multline*}
   (a_i)_{i\in \{1,\dots,m\}},\ a_i \in \{Radj,\ Ladj,\ Ru,\ Rco,\ Lu,\ Lco\},\\ a_m(\dots a_2(a_1(f))\dots)\text{ exists if and only if }a_m(\dots a_2(a_1(g))\dots)\text{ exists,} 
\end{multline*} and this for any choice of adjoints, units and counits.
\end{mydef}
We will show that dualizability properties are preserved by isomorphisms and ``higher mating" defined in Definition \ref{def:mates}. Let us describe formally this second notion.
\begin{prop}
Let $f:x\to y$ be a $k$-morphism in an $n$-category $\CC$ with adjoint $(f^R, \varepsilon, \eta)$. Then for any other $k$-morphisms $g: z \to x$ and $h: z \to y$, one has an equivalence between the $(n-k-1)$-categories of $(k+1)$-morphisms: $$\Phi^f_{g,h}: \left\{ \begin{array}{ccc}\Hom_\CC(f \circ_k g, h) &\tilde\to& \Hom_\CC(g, f^R \circ_k h)\\
N: f \circ_k g \to h &\mapsto& (\Id_{f^R}\circ_k N) \circ_{k+1} (\eta \circ_k \Id_g)\\
\text{$k+j$-morphism } \alpha &\mapsto& (\Id^j_{f^R}\circ_k \alpha) \circ_{k+1} (\Id^{j-1}_\eta \circ_k \Id^j_g)
\end{array}\right\}$$
Similarly, for any $g:x \to z$ and $h:y \to z$, one gets an equivalence: $$\Psi^f_{g,h}= (-\circ_k \Id_{f})\circ_{k+1} (\Id_g\circ_k\eta):\Hom_\CC(g \circ_k f^R, h) \tilde\to \Hom_\CC(g, h \circ_k f)$$
\end{prop}
\begin{myproof} Its inverse is given by:
$$(\Phi^f_{g,h})^{-1}: \left\{ \begin{array}{ccc}\Hom_\CC(g, f^R \circ_k h)&\tilde\to& \Hom_\CC(f \circ_k g, h) \\
M: g \to f^R \circ_k h &\mapsto& (\varepsilon\circ_k \Id_h) \circ_{k+1} (\Id_f \circ_k M) \\
\text{$(k+j)$-morphism } \beta &\mapsto& (\Id_\varepsilon^{j-1}\circ_k \Id^j_h) \circ_{k+1} (\Id^j_f \circ_k \beta)
\end{array}\right\}$$
The composition $\Phi^f_{g,h} \circ (\Phi^f_{g,h})^{-1}$ (resp. $(\Phi^f_{g,h})^{-1}\circ \Phi^f_{g,h}$) is post- (resp. pre-) composition by a snake identity. Similarly, $(\Psi^f_{g,h})^{-1} = (\Id_g \circ_k \varepsilon)\circ_{k+1} (- \circ_k \Id_{f^R})$.
\end{myproof}
\begin{mydef}\label{def:mates}
For a $(k+1)$-morphism $N: f \circ_k g \to h$, we say that $N$ and $\Phi^f_{g,h}(N)$ are \emph{mates}. For a higher morphism $\alpha$ in $\Hom_\CC(f \circ_k g, h)$, we say that $\alpha$ and $\Phi^f_{g,h}(\alpha)$ are \emph{higher mates}. Similarly, for $N,\ \alpha$ in $\Hom_\CC(g \circ_k f^R, h)$ we call $N$ and $\Psi^f_{g,h}(N)$ mates, and $\alpha$ and $\Psi^f_{g,h}(\alpha)$ higher mates. More generally we say that $N$ and $M$ are mates (resp. $\alpha$ and $\beta$ are higher mates) if they can be linked by a chain of matings (resp. higher matings) and isomorphisms.

For a $k$-morphism $f$, we say that $g$ is \emph{obtained from $f$ by whiskering} if it can be written as a composition of $f$ with identities of lower morphisms.
Note that if $\alpha$ and $\beta$ are higher mates, their are both obtained from the other by whiskering.
\end{mydef}
\begin{prop}\label{prop:coherenceDualiz} Let $f$ and $g$ be $k$-morphisms in an $n$-category. Then:
\begin{enumerate}
\item $f \simdual f$.
\item If $f \overset{\varphi}{\simeq} g$ are isomorphic, then $f \simdual g$.
\item If $f = g \circ_k h$ for an isomorphism $h$, then $f \simdual g$.
\item If $f$ and $g$ are higher mates, then $f \simdual g$.\end{enumerate}
\end{prop}
\begin{myproof} What we have to prove for point 1 is that existence of higher adjoints in the adjunctibility data does not depend on the choices made in the adjunctions. This kind of results are known as coherence statement in the literature and is usually stated as contractibility of a space of dualizability data, see \cite[Lemma 4.6.1.10]{LurieHA} or \cite{RiehlVerityCoherence}, though our statement here is more elementary.

Let us discuss only right adjoints and right units below, every other notion being related by taking appropriate opposite categories.

The adjunctibility data of $f$ is unique up to isomorphism. This means that if $g_1$ and $g_2$ are both right adjoints to $f$, with units $u_1$ 
and $u_2$, then there is an isomorphism $\varphi: g_1 \tilde\to g_2$ such that $u_2$ is isomorphic to $(\varphi \circ_k \Id_f)\circ_{k+1} u_1$.

We observe that any choice of right adjoint or right unit of $f$ are related by a sequence of either point 2 (isomorphism) or point 3 (composing with an isomorphism). The strategy of the proof is to show more generally that if $f$ and $g$ are related by a sequence of either point 2 or point 3, then $f$ is right adjunctible if and only $g$ is, and their adjoints, units and counits are again related by a sequence of either point 2 or point 3. The result then follows by induction on $m$ the number of letters in $a_m(\dots a_2(a_1(f))\dots)$. For point 4 we show similarly that if $f$ is related to $g$ by a sequence of point 2, 3 and 4, then so are their adjoints, units and counits.

\begin{enumerate}\setcounter{enumi}{1}
\item If $f \overset{\varphi}{\simeq} g$ are isomorphic, then $f$ has a right (resp. left) adjoint if and only if $g$ does, in which case one can choose $Radj(g) = Radj(f)$, $Ru(g) = (\Id_{Radj(f)} \circ_k \varphi) \circ_{k+1} Ru(f)$ and $Rco(g) = (\varphi^{-1}\circ_k \Id_{Radj(f)}) \circ_{k+1} Rco(f)$.
\item If $f = g \circ_k h$ is obtained as a composition, then $f$ has a right (resp. left) adjoint as soon as $g$ and $h$ do, in which case one can choose $Radj(f) = Radj(h) \circ_k Radj(g)$, $Ru(f) = (\Id_{Radj(h)} \circ_k Ru(g) \circ_k \Id_h) \circ_{k+1} Ru(h)$ and $Rco(f) = (\Id_g \circ_k Rco(h) \circ_k \Id_{Radj(g)}) \circ_{k+1} Rco(h)$. In particular, if $h$ is an isomorphism, then $g \simeq f \circ_k h^{-1}$, and $f$ has a right (resp. left) adjoint if and only if $g$ does.
\item If $f = g \circ_j h$ is obtained as a composition in the direction of $j$-morphisms for $j<k$, then $f$ has a right (resp. left) adjoint as soon as $g$ and $h$ do, in which case one can choose $Radj(f) = Radj(g) \circ_j Radj(h)$, $Ru(f) = Ru(g) \circ_j Ru(h)$ and $Rco(f) = Rco(g)\circ_j Rco(h)$. In particular, if $h$ is an identity of a lower morphism, then $f$ has a right (resp. left) adjoint as soon as $g$ does. So, if $f$ and $g$ are higher mates, they both can be obtained as composition of the other with identities of lower morphisms, and $f$ has a right (resp. left) adjoint if and only if $g$ does.
\end{enumerate}
\end{myproof}
We can now describe the redundancy in the dualizability data:
\begin{prop} \label{propRedundancyDD}
 Let $f$ be a $k$-morphism in an $n$-category $\CC$, suppose that $Radj(f),\ Radj(Rco(f))$ and $Radj(Ru(f))$ exist, then:
 \begin{enumerate}
     \item $f$ is 1-adjunctible, and one can choose $Ladj(f) = Radj(f)$, $Lu(f) = Radj(Rco(f))$ and $Lco(f)= Radj(Ru(f))$.
     \item $Rco(Ru(f)) \simdual Ru(Rco(f))$.
 \end{enumerate}
 Suppose moreover than $Radj(Lco(f))$ and $Radj(Lu(f))$ exist, then:
 \begin{enumerate}\setcounter{enumi}{2}
     \item $f$ is 2-adjunctible, and $Rco(f) \simdual Radj(Radj(Rco(f)))$ and $Ru(f) \simdual Radj(Radj(Ru(f)))$.
 \end{enumerate}
 In particular if $f=X$ is an object in a symmetric monoidal $n$-category, then:
  \begin{enumerate}\setcounter{enumi}{3}
     \item $X$ is 1-adjunctible if and only if it has a dual. It is 2-adjunctible if and only if $ev_X := Rco(X)$ and $coev_X := Ru(X)$ have right adjoints. More generally, it is $m$-adjunctible if and only if $Radj(Rco^k(Ru^{m-1-k}(X)))$ exist for all $0 \leq k \leq m-1$.
 \end{enumerate}
\end{prop}
\begin{myproof}
Point 1 is \cite[Remark 3.4.22]{LurieCob}, or \cite[Lemma 20.1]{SP}. One directly checks that the right adjoints of the right counit and unit satisfy the snake relations, because taking right adjoints behaves well with composition, and exhibit $Radj(f)$ as the left adjoint of $f$. 

Point 2 is \cite[Proposition 3.4.21]{LurieCob}. It is shown that $Rco(Ru(f))$ and $Ru(Rco(f))$ are higher mates, so in particular $Rco(Ru(f)) \simdual Ru(Rco(f))$.

Point 3 is \cite[Lemma 7.11]{JFS}. 
Let us recall the argument to illustrate what we mean by redundancy in the dualizability data. First by point 1 we can take $Lco(f)$ and $Lu(f)$ to be the right adjoints of $Ru(f)$ and $Rco(f)$. Now by a left-handed version of point 1 we can also obtain a unit and counit for the right adjunction of $f$ as the right adjoints of $Lco(f)$ and $Lu(f)$. 
Therefore both $Rco(f)$ and its double right adjoint are counits for the right adjunction of $f$: there is a redundancy in the dualizability data. By the proposition above, they have the same dualizability properties, and similarly for the right unit.

Point 4: by \cite[Corollary 7.12]{JFS} we have to check that $X$ is dualizable and that its evaluation and coevaluation maps are $(m-1)$-times right adjunctible, i.e. that $Radj(a_1(\dots a_{m-1}(X)\dots)$ exist for any $(a_i)_i\in \{Rco, Ru\}^{m-1}$. Using point 2, we know that $Ru$ and $Rco$ commute as far as existence of adjoints is concerned, so there are only $m$ different $m-1$-morphisms whose adjunctibility should be checked, $Rco^k(Ru^{m-1-k}(X))$, $0 \leq k \leq m-1$. 
\end{myproof}

\subsubsection{Oplax dualizability data}\label{subsectOplax} We investigate the proof of Theorem \ref{thmCaracOplaxDualiz} and explain how to get from adjunctibility data in $\CC$ to dualizability data in $\CC^\to$.
\begin{myth}[Johnson-Freyd--Scheimbauer] \label{thmDDinOplax} Let $f= (s_f, t_f, f^\#) : a = (s_a, t_a, a^\#) \to b = (s_b,t_b,b^\#)$ be a $k$-morphism in $\CC^\to$ so $s_f: s_a \to s_b$ and $t_f: t_a \to t_b$ are $k$-morphism in $\CC$, and $f^\#$ is a $k+1$-morphism in $\CC$ from the composition of some whiskerings of $b^\#$ and $s_f$ to the composition of some whiskerings of $t_f$ and $a^\#$. Then:
\begin{center}
$f$ has a right adjoint in $\CC^\to$ if and only if $s_f,\ t_f$ and $f^\#$ have right adjoints in $\CC$.
\end{center}
 In this case: \begin{enumerate}[\footnotesize $\bullet$]
     \item $Radj(f) = (Radj(s_f),Radj(t_f), g)$ where $g$ is a mate of $Radj(f^\#)$,
     \item $Ru(f) = (Ru(s_f), Ru(t_f), u)$ where $u$ is a higher mate of $Rco(f^\#)$, and 
     \item $Rco(f) = (Rco(s_f), Rco(t_f), v)$ where $v$ is a higher mate of $Ru(f^\#)$.
 \end{enumerate}
In particular, if we only look at the right dualizability data, and only take right adjoints once, then:
\begin{center}
$\forall i,j\in\Nn,\ Radj(Rco^i(Ru^j(f)))$ exists if and only if $Radj(Rco^i(Ru^j(s_f)))$, $Radj(Rco^i(Ru^j(t_f)))$ and $Radj(Ru^i(Rco^j(f^\#)))$ exist.
\end{center}
\end{myth}
\begin{myproof}
The description of the right adjunctibility of a morphism in $\CC^\to$ is \cite[Proposition 7.13]{JFS}, in the oplax case.

For the last statement, remember that $u\simdual Rco(f^\#)$ and $v\simdual Ru(f^\#)$ by Proposition \ref{prop:coherenceDualiz}.4. The first statement implies by induction that $Rco^i(Ru^j(f))$ is of the form $(s,t,w)$ where $s \simdual Rco^i(Ru^j(s_f))$, $t\simdual Rco^i(Ru^j(t_f))$ and $w\simdual Ru^i(Rco^j(f))$.

Indeed increasing $j$ for $i=0$ we have $Ru^{j+1}(f) = (Ru(s),Ru(t), U)$ which have the same dualizability properties as respectively $(Ru^{j+1}(s_f),Ru^{j+1}(t_f),Rco^{j+1}(f^\#))$. Increasing $i$ we have $Rco^{i+1}Ru^j(f) = (Rco(s), Rco(t), V)$ which have the same dualizability properties as respectively $(Rco^{i+1}Ru^j(s_f), Rco^{i+1}Ru^j(t_f), Ru^{i+1}Rco^j(f^\#))$.
\end{myproof}
\begin{myex}[$k=0$]\label{exmplDualinCto} An object $f=(X,Y,A:X\to Y)$ of $\CC^\to$ is dualizable if and only if $X$ and $Y$ are dualizable, and $A$ has a right adjoint $Radj(A)$. Then:
\begin{enumerate}[\footnotesize $\bullet$]
     \item $f^* = (X^*, Y^*, Radj(A)^* := (\Id_{Y^*}\otimes\, ev_X) \circ (\Id_Y^* \otimes Radj(A) \otimes \Id_{X^*})\circ (coev_Y \otimes \Id_{X^*}))$,
     \item $coev_f = (coev_X, coev_Y, (Rco(A)\otimes \Id_{\Id_{Y^*}})\circ_1 \Id_{coev_Y})$, and 
     \item $ev_f = (ev_X, ev_Y, \Id_{ev_X}\circ_1(Ru(A) \otimes \Id_{\Id_{X^*}}))$.
 \end{enumerate} 
A surprising consequence of this result is that if $f$ is 2-dualizable, the right counit and unit of $A$ are biadjoints up to isomorphisms and mating. A drawing for this is given in Figure \ref{unitCounitCompose}.
\end{myex}

\subsection{Cobordism Hypotheses}
The Cobordism Hypothesis describes fully extended Topological Quantum Field Theories with values in a higher category $\CC$ in terms of fully dualizable objects of $\CC$. We also recall relative versions that describes relative TQFTs, and a non-compact version that describes partially-defined TQFTs. The Cobordism Hypothesis was formulated in \cite{BaezDolan}. A sketch of proof was given in \cite{LurieCob}, a more formal version is work in progress of Schommer-Pries. An independent proof of a more general result appears in the preprint \cite{GradyPavlov}. Another independent proof using factorization homology is work in progress, see \cite{AyalaFrancisCH}.
\begin{mythq}[The Cobordism Hypothesis, Thm 2.4.6 and 2.4.26 in \cite{LurieCob}]
Let $\CC$ be a symmetric monoidal $n$-category. Evaluation at the point induces equivalences of $\infty$-groupoids
$$\Fun^\otimes(\Bord_n^{fr},\CC)\simeq (\CC^{fd})^\sim$$
between framed fully extended $n$-TQFTs with values in $\CC$ and the underlying $\infty$-groupoid of the subcategory of fully dualizable objects of $\CC$, and 
$$\Fun^\otimes(\Bord_n,\CC)\simeq ((\CC^{fd})^\sim)^{SO(n)}$$
between oriented fully extended $n$-TQFTs with values in $\CC$ and $SO(n)$-homotopy-fixed-points in $(\CC^{fd})^\sim$, where $SO(n)$ acts on the $n$-category $\Bord_n^{fr}$ by changing the framing and therefore on $(\CC^{fd})^\sim$ by the first equivalence.
\end{mythq}
For $X \in \CC$ a fully dualizable object, we denote $\ZZ_X$ (a choice of representative of) the associated fully extended framed $n$-TQFT.

\subsubsection{The relative Cobordism Hypothesis}
Lurie proposes a result classifying his notion of domain wall.
\begin{mythq}[Theorem 4.3.11 and Example 4.3.23 in \cite{LurieCob}] Let $\CC$ be a symmetric monoidal $n$-category and $X,Y \in \CC^{fd}$. There is a bijection between isomorphism classes of framed domain walls between $\ZZ_X$ and $\ZZ_Y$ and isomorphism classes of fully dualizable 1-morphisms $f:X\to Y$, given by evaluation at the interval with an interface point in the middle.
\end{mythq}
There is an oriented version asking that $f$ preserves orientation structures.

On the other hand, \cite{JFS}'s notions of a twisted quantum field theory are already classified by the usual Cobordism Hypothesis. Note however that \cite[Definition 5.16]{JFS} is surprisingly strict because it demands that the source and target of the functor $\RR:\Bord_{n-1}\to\CC$ agree strictly with $\ZZ_1$ and $\ZZ_2$. Equivalently, we could have asked that $\RR$ comes equipped with isomorphisms $s(\RR)\simeq \ZZ_1$ and $t(\RR)\simeq \ZZ_2$. In both cases, it is clear that the Cobordism Hypothesis does not apply on the nose. The fix is easy.
\begin{mydef}
Let $\CC$ be a symmetric monoidal $n$-category and $X,Y \in \CC$. Denote $(\CC^\to)^\sim_{X,Y}$ the homotopy pullback $$\begin{tikzcd}
    (\CC^\to)^\sim_{X,Y} \ar[r]\ar[d, very near start, right = 10pt, "\quad \lrcorner_h"] & (\CC^\to)^\sim \ar[d,"{s,t}"]\\
    \ast \ar[r, "{X,Y}"] & (\CC^\sim)^{\times 2}
\end{tikzcd}.$$
Similarly, for $\ZZ_1,\ZZ_2: \Bord_{n-1}^{fr}\to\CC$ denote $\Fun^\otimes(\Bord_{n-1}^{fr},\CC^\to)_{\ZZ_1,\ZZ_2}$ the homotopy pullback $$\begin{tikzcd}
    \Fun^\otimes(\Bord_{n-1}^{fr},\CC^\to)_{\ZZ_1,\ZZ_2} \ar[r]\ar[d, very near start, right = 10pt, "\quad \lrcorner_h"] & \Fun^\otimes(\Bord_{n-1}^{fr},\CC^\to) \ar[d,"{s,t}"]\\
    \ast \ar[r, "{\ZZ_1,\ZZ_2}"] & (\Fun^\otimes(\Bord_{n-1}^{fr},\CC))^{\times 2}
\end{tikzcd}$$ called the space of framed oplax-$\ZZ_1$-$\ZZ_2$-twisted-$(n-1)$-TQFTs.

Note that both are also strict pullbacks as taking source and target induces a fibration of spaces.
\end{mydef} 
\begin{corr}[of the Cobordism Hypothesis]\label{CHforOplax}
Let $\CC$ be a symmetric monoidal $n$-category and $X,Y \in \CC$. Choose $\ZZ_X,\ZZ_Y: \Bord_{n-1}^{fr}\to\CC$ two TQFTs associated with $X$ and $Y$ by the cobordism hypothesis. Evaluation at the point induces an equivalence
$$ \Fun^\otimes(\Bord_{n-1}^{fr},\CC^\to)_{\ZZ_X,\ZZ_Y} \simeq (\CC^\to)^\sim_{X,Y}\ .$$
\end{corr}
\begin{myproof}
    The cobordism hypothesis on $\CC$ and $\CC^\to$ gives a commutative diagram of horizontal equivalences 
    $$\begin{tikzcd}
        \Fun^\otimes(\Bord_{n-1}^{fr},\CC^\to)\ar[d,"{s,t}"] \ar[r, "ev_{pt}"] & (\CC^\to)^\sim \ar[d,"{s,t}"]\\
        (\Fun^\otimes(\Bord_{n-1}^{fr},\CC))^{\times 2}\ar[r, "ev_{pt}\times ev_{pt}"] & (\CC^\sim)^{\times 2}\\
        \ast\ar[u,"{\ZZ_X,\ZZ_Y}"] \ar[r] & \ast \ar[u,"{X,Y}"]
    \end{tikzcd}$$ inducing an equivalence between homotopy pullbacks.
\end{myproof}
\begin{rmk}\label{rmkOrientedTwistedCH}
There is an oriented version as well. The maps $s,t: \Fun^\otimes(\Bord_{n-1}^{fr},\CC^\to) \to \Fun^\otimes(\Bord_{n-1}^{fr},\CC)$ are $SO(n-1)$-equivariant because $SO(n-1)$ acts on the source $\Bord_{n-1}^{fr}$. Therefore the maps $s,t: (\CC^\to)^\sim \to \CC^\sim$ are also equivariant, and descend to maps between the $SO(n-1)$-homotopy-fixed-points $s,t: (\CC^\to)^{\sim, SO(n-1)} \to \CC^{\sim, SO(n-1)}$. Given two objects $X,Y\in\CC$ equipped with $SO(n-1)$-homotopy-fixed point structure, one can reproduce exactly the whole paragraph above and define $(\CC^\to)^{\sim, SO(n-1)}_{X,Y}$ as a pullback. We get $$ \Fun^\otimes(\Bord_{n-1},\CC^\to)_{\ZZ_X,\ZZ_Y} \simeq (\CC^\to)^{\sim, SO(n-1)}_{X,Y}$$ by the same proof, using the oriented cobordism hypothesis.    
\end{rmk} 
\begin{rmk} \label{rmkWillStewart}
Results to-appear in the PhD of William Stewart show that if we assume that the source and target objects $X$ and $Y$ are fully dualizable then a morphism $f:X\to Y$ is $(n-1)$-oplax dualizable if and only if it is $(n-1)$-dualizable. In particular, if we restrict the notion of oplax twisted TQFTs to the case where the ``twisting" theories $\ZZ_1$ and $\ZZ_2$ extend to $\Bord_n$, which is the setting in \cite{FreedTeleman}, then this notion, using the cobordism hypothesis twice, is equivalent to Lurie's notion of domain walls.
\end{rmk}

\subsubsection{Non-compact TQFTs}\label{subsectNonCompact} To study non-semisimple variants of Witten--Reshetikhin--Turaev TQFTs, we will be interested in theories defined on a restricted class of cobordisms, namely where top-dimensional cobordisms have non-empty outgoing boundary in every connected component.

Lurie's sketch of proof of the cobordism hypothesis is done by induction on the handle indices allowed. One starts with only opening balls, then allows more and more complex cobordisms. Eventually one allows every cobordisms but closing balls, namely cobordisms with outgoing boundary in every connected component. Finally one allows every cobordism, and obtain a TQFT. We call it a non-compact TQFT when we stop at this ante-last step. Lurie's proof then gives an algebraic criterion classifying these. 

We follow \cite[Section 3.4]{LurieCob} and state the results there in a form fitted for our use. It should be noted that the proofs of the statements below are not very formal.
\begin{mydef}
Let $\Bord_n^{fr,nc} \subseteq \Bord_n^{fr}$ denote the subcategory where $n$-dimensional bordisms have non-empty outgoing boundary in every connected component.

A \emph{framed fully extended non-compact $n$-TQFT} with values in a symmetric monoidal $n$-category $\CC$ is a symmetric monoidal functor $\ZZ:\Bord_n^{fr,nc}\to \CC$. 
\end{mydef}
Lurie defines in \cite[Definition 3.4.9]{LurieCob} an $n$-category $\FF_k$ of $\leq n$-dimensional bordisms where all $n$-manifolds are equipped with a decomposition into handles of index $\leq k$. Here bordisms are actually equipped with a framed function without certain kinds of critical points.

We denote $\alpha_k^m = D^k\times D^{m-k}: S^{k-1} \times D^{m-k} \to D^k \times S^{m-k-1}$ the $m$-dimensional index $k$ handle attachment, seen as an $m$-morphism in $\Bord_m^{fr}$, or in $\FF_k$ if $m=n$. Let $x = S^{k-2}\times D^{n-k}$, $y=D^{k-1}\times S^{n-k-1}$ seen as $(n-2)$-morphisms $\emptyset \to S^{k-2}\times S^{n-k-1}$ in $\Bord_{n-1}^{fr}$. Note that for $1\leq k\leq n$, $\alpha_{k-1}^{n-1}: x \to y$ and $\alpha_{n-k}^{n-1}:y\to x$. Then, $\alpha_{k-1}^n$ can be seen (up to higher mating) as a morphism $\Id_x \to \alpha_{n-k}^{n-1}\circ \alpha_{k-1}^{n-1}$ and $\alpha_k^n$ as a morphism $\alpha_{k-1}^{n-1} \circ \alpha_{n-k}^{n-1} \to \Id_y$, and they form a unit/counit pair in $\FF_k$, see \cite[Claim 3.4.17]{LurieCob}. Namely, $Radj(\alpha_{k-1}^{n-1})=\alpha_{n-k}^{n-1}$, $Ru(\alpha_{k-1}^{n-1})=\alpha_{k-1}^n$ and $Rco(\alpha_{k-1}^{n-1})=\alpha_k^n$, or in our diagrammatic notation:
\begin{center}
    \begin{tikzpicture}[baseline = 10pt, xscale = 1.3, yscale = 1.3]
\node (A) at (1,1) {$\alpha_{k-1}^{n-1}$};
\node (R) at (2,0) { $\alpha_k^n\ \vert\ \alpha_{k-1}^n$};
\draw[->] (A) -- (R) node[near end,above right]{$\alpha_{n-k}^{n-1}$};
\end{tikzpicture}\ .
\end{center}
 By induction, $Rco^{k}(Ru^{m-k}(pt))=\alpha_k^m$.
\begin{mythq}[Index-$k$ cobordism hypothesis, Lurie]
A symmetric monoidal functor $\ZZ_0: \Bord_{n-1}^{fr} \to \CC$ extends to $\ZZ: \FF_k \to \CC$, $1\leq k\leq n$, if and only if the images of every $(n-1)$-dimensional handle of index $\leq k-1$ is right adjunctible. 

This extension is essentially unique: there is an isomorphism $\ZZ\Rightarrow\ZZ'$ between any two such extensions. This isomorphism may not be the identity on $\Bord_{n-1}^{fr}$.
\end{mythq}
\begin{myproof}[Sketch]
For $k=0$, one can extend $\ZZ_0: \Bord_{n-1}^{fr} \to \CC$ with any $n$-morphism $\ZZ(\alpha_0^n):1 \to \ZZ_0(S^{n-1})$, see \cite[Claim 3.4.13]{LurieCob}. Note that Lurie works in the unoriented case there, and demands on $O(n)$-equivariant morphism, and we look at the framed case.

Now, for $1\leq k\leq n$, a symmetric monoidal functor $\ZZ_0: \FF_{k-1} \to \CC$ extends to $\ZZ: \FF_k \to \CC$ if and only if $\alpha_{k-1}^n$ is mapped to a unit of an adjunction between $\alpha_{k-1}^{n-1}$ and $\alpha_{n-k}^{n-1}$, see \cite[Proposition 3.4.19]{LurieCob}. In this case, the extension is essentially unique, and $\alpha_k^n$ is mapped to the counit of the adjunction.

For $k=1$, this gives little choice for the $n$-morphism $\ZZ(\alpha_0^n)$: it has to be the unit of an adjunction and is therefore determined up to isomorphism. Then, $\alpha_1^n$ will be sent to the counit. 

For $k\geq 2$, we want $\ZZ(\alpha_{k-1}^n)$, which is so far defined as the counit of the adjunction between $\ZZ(\alpha_{k-2}^{n-1})$ and $\ZZ(\alpha_{n-k+1}^{n-1})$, to be also the unit of the adjunction between $\ZZ(\alpha_{k-1}^{n-1})$ and $\ZZ(\alpha_{n-k}^{n-1})$. This in particular implies that the $(n-1)$-dimensional handle of index $k-1$ is right adjunctible, as stated in the conjecture. For the converse, we use \cite[Proposition 3.4.20]{LurieCob} (which we recalled in Proposition \ref{propRedundancyDD}.2) which states that provided the adjunction exists, $\alpha_{k-1}^n$ must map to (a higher mate of) the unit.
\end{myproof}
\begin{mydef}
Let $\CC$ be a symmetric monoidal $n$-category. An object $X$ in $\CC$ is said $(n,k)$-dualizable if it is $n-1$-dualizable and the $k$ following $n-1$-morphisms $Ru^{n-1}(X)$, $Rco(Ru^{n-2}(X))$, $\dots$, $Rco^{k-1}(Ru^{n-k}(X))$ have right adjoints. We say $X$ is non-compact-$n$-dualizable if it is $(n,n-1)$-dualizable.

For example, for $n=3$, $k=2$, we want $X$ to have a dual $(X^*, ev_X, coev_X)$, both its evaluation and coevaluation maps to have right adjoints $(ev_X^R, a, b)$ and $(coev_X^R, c, d)$, and the unit and counit of the right adjunction of the coevaluation to have right adjoints $c^R$ and $d^R$.
\end{mydef}
We can now state the non-compact version of the cobordism hypothesis, which we will assume in Section \ref{sectAppli}. A formal proof is work-in-progress of Schommer-Pries.
\begin{mythq}[Non-compact Cobordism Hypothesis]
Let $\CC$ be a symmetric monoidal $n$-category, $n\geq 2$. There is a bijection between isomorphism classes of framed fully extended non-compact $n$-TQFTs with values in $\CC$ and isomorphism classes of non-compact-$n$-dualizable objects of $\CC$, given by evaluation at the point.
\end{mythq}
There is an oriented version as well, stating that oriented non-compact theories are classified by $SO(n)$-homotopy fixed points in the space of non-compact-$n$-dualizable objects.
\section{Dualizability of the unit inclusion} \label{sectDD} 
Let $\VV\in \BrTens$ be a braided tensor category. We consider the inclusion of the unit $\eta: \Vect_\FK \to\VV$. It is a braided monoidal functor and we define an associated $\Vect_\FK$-$\VV$-central algebra $\AA_\eta$, which is simply the category $\VV$ seen as the regular right $\VV$-module, see Definition \ref{defAF}. We study the dualizability of this 1-morphism in $\BrTens$. First, we recall some context and develop some properties of bimodules induced by functors. Then we describe all the dualizability data explicitly and give criteria for dualizability.

\subsection{Cocomplete braided tensor categories}\label{subsectContext}
We will work in the even higher Morita category $\Alg_2(\Pr)$. This category have been formally defined in \cite{JFS}, and described more explicitly in \cite{BJS} under the name $\BrTens$. 
\subsubsection{Cocomplete categories} 
We begin by recalling some properties of the 2-category $\Pr$. Let $\FK$ be a field of characteristic zero. 

 Let $\Cat_\FK$ denote the 2-category of small $\FK$-linear categories, and $\Pr$ denote the 2-category of cocomplete locally presentable $\FK$-linear categories \cite[Defintion 1.17]{AdamekRosicky}, $\FK$-linear cocontinuous functors and $\FK$-linear natural transformations, equipped with the Kelly tensor product $\boxtimes$. We denote $\Free = \Hom_{\Cat_\FK}((-)^{op},\Vect_\FK) : \Cat_\FK \to \Pr$ the symmetric monoidal free cocompletion functor. Its essential image is denoted $\Bimod_\FK$\footnote{The name comes from the Eilenberg--Watts theorem which describes cocontinuous functors between categories of modules over two algebras as bimodules.}.
 
 An object $C\in\CC$ is called compact-projective (which we abbreviate cp) if the functor $\Hom_\CC(C,-)$ is cocontinuous. The category $\CC$ is said to have enough compact-projectives if its full subcategory $\CC^{cp}$ of cp objects generates $\CC$ under colimits, or equivalently if the canonical functor $\Free(\CC^{cp})\to\CC$ is an equivalence, or if it lies in $\Bimod_\FK$. A monoidal category $\CC$ is called cp-rigid if it has enough cp and all its cp objects are left and right dualizable.
\begin{prop}
A 1-morphism $F:\CC\to \DD$ in $\Pr$ between two categories with enough cp has a cocontinuous right adjoint if and only if it preserves cp.
\end{prop}
\begin{myproof}
If $F^R$ is cocontinuous then for $C\in\CC^{cp}$ and $D=\colim_i D_i$ obtained as a colimit, $$\begin{array}{ccl}
\Hom_\DD(F(C),D)&\simeq&\Hom_\CC(C,F^R(D))\\
&\overset{F^R\text{ cocont}}{\simeq}& \Hom_\CC(C,\colim_i F^R(D_i)) \\
&\overset{C\text{ cp}}{\simeq}& 
\colim_i \Hom_\CC(C,F^R(D_i))\simeq\colim_i \Hom_\DD(F(C),D_i)
\end{array}$$ and $F(C)$ is compact-projective.

The other direction is a classical construction, see \mbox{\cite[Lemma 2.10]{BDSPV}}.
\end{myproof}
\begin{rmk}
The condition of $F$ preserving cp, namely of $F$ being the cocontinuous extension of a functor $f$ on the subcategories of cp objects, is very similar to that of a bimodule being induced by a functor in Section \ref{sectFunBimod}. When $F$ preserves cp then $F^R$ is associated with the ``mirrored" bimodule induced by $f$, with unit induced by $f$, and counit induced by composition in $\DD$. This is again very similar to what happens in Section \ref{sectFunBimod}.
\end{rmk}

\subsubsection{$\BrTens$ and \sc{Alg}${}_2(\Pr)$} \label{subsBrTensAlg2}
The higher Morita $n$-category $\Alg_n(\mathcal{S})$ associated with an $\infty$-category $\mathcal{S}$ was introduced in \cite{Haugseng} using a combinatorial/operadic description. A pointed version was introduced in \cite{Scheimbauer} using very geometric means, namely factorization algebras. This geometric description allows for a good description of dualizability in $\Alg_n(\mathcal{S})$, but the pointing prevents any higher dualizability, see \cite{GS}. Even higher Morita categories are defined in \cite{JFS}, for pointed and unpointed versions. They form an $(n+k)$-category $\Alg_n(\mathcal{S})$ for $\mathcal{S}$ a $k$-category. We consider the unpointed even higher Morita 4-category $\Alg_2(\Pr)$, which we denote $\BrTens$ for reasons that will be made explicit below. Even though we are not formally in this context, we will sometimes use factorization algebra drawings to illustrate our point.

One represents an $\mathbb E_2$-algebra $\VV$ as \VVrpz\ and the $\Vect_\FK$--$\VV$-algebra $\AA_\eta$ as \AArpz.
Let us recall the description of $\BrTens$ from \cite{BJS}.
\begin{mydef}[Section 2.4 in \cite{BJS}]
An object $\VV$ of $\BrTens$ is a locally presentable cocomplete $\FK$-linear braided monoidal category. We call these \emph{braided tensor categories} here, even though this name has many uses. Equivalently, it is an $\mathbb E_2$-algebra in $\Pr$.
\end{mydef}
In the factorization algebra picture \begin{tikzpicture}[baseline = -20pt, yscale=-1]
\fill[blue!30] (0,0.5) rectangle (1,1);
\fill[red!15] (0,0) rectangle (1,0.5);
\draw[very thick] (0,0.5)--(1,0.5);
\end{tikzpicture}, one can read that a 1-morphism between two braided tensor categories $\VV$ (in red, on top) and $\WW$ (in blue, below) is a monoidal category $\AA\in \Pr$ (the horizontal line) with a top $\VV$-action and a bottom $\WW$-action that commute with respect to each other and that commute with the monoidal structure of $\AA$ in a coherent way. Note that as $\AA$ is monoidal, such an $\VV$-action $\rhd$ is determined by a monoidal functor $\VV\to \AA$ that maps $V$ to $V \rhd \bm1_\AA$. See \cite[Figure 2]{BJS}.
\begin{mydef}[Definition-Proposition 3.2 in \cite{BJS}]
A 1-morphism between $\VV$ and $\WW$ in $\BrTens$ is a \emph{$\VV$-$\WW$-central algebra} $\AA$. Namely, it is an monoidal category $\AA \in \Pr$ equipped with a braided monoidal functor $$(F_\AA,\sigma^\AA): \VV \boxtimes \WW^{\sigma op} \to Z(\AA)$$ to the Drinfeld center of $\AA$.

Remember that the Drinfeld center of $\AA$ has objects pairs $(y,\beta)$ where $y$ is an object of $\AA$ and $\beta: -\otimes y\ \tilde{\Rightarrow}\ y\otimes -$ is a natural isomorphism. Here $F_\AA$ gives the object and $\sigma^\AA$ gives the half braiding. We denote $V\rhd A := F_\AA(V) \otimes A$ and $A \lhd V := A \otimes F_\AA(V)$.
\end{mydef} 
Composition of 1-morphism is relative tensor product over the corresponding braided tensor category, see \cite[Section 3.4]{BJS}

Again in the factorization algebra picture \begin{tikzpicture}[baseline = -20pt, yscale=-1]
\fill[blue!30] (0,0.5) rectangle (1,1);
\fill[red!15] (0,0) rectangle (1,0.5);
\draw[very thick] (0,0.5)--(0.5,0.5);
\draw[densely dashed, very thick] (0.5,0.5)--(1,0.5);
\node at (0.5,0.5) {$\bullet$};
\end{tikzpicture}, a 2-morphism $\MM$ between two $\VV$-$\WW$-central algebras $\AA$ and $\BB$ is a $\AA$-$\BB$-bimodule category where $\VV$ (resp. $\WW$) acts similarly when acting through $\AA$ or through $\BB$.
\begin{mydef}[Definition 3.9 in \cite{BJS}]
A 2-morphism between $\AA$ and $\BB$ in $\BrTens$ is a \emph{$\VV$-$\WW$-centered $\AA$-$\BB$-bimodule category}. Namely, it is an $\AA$-$\BB$-bimodule category $\MM$ equipped with natural isomorphisms $$\eta_{v,m}: F_\AA(v)\rhd m \ \tilde{\to}\ m\lhd F_\BB(v)\ , \quad v\in\VV\, ,\ m\in\MM\ ,$$ satisfying coherences with tensor product in $\VV$ and with the half braidings in $\AA$ and $\BB$.
\end{mydef}
Horizontal and vertical composition are again relative tensor product over the corresponding monoidal category.
\begin{mydef}[Section 3.6 in \cite{BJS}]
A 3-morphism $F:\MM\to\NN$ is a functor of $\AA$-$\BB$-bimodules categories that preserves the $\VV$-$\WW$-centered structure.

A 4-morphism $\eta: F \Rightarrow G$ is a natural transformation of bimodule functors.
\end{mydef}

\subsubsection{Previous dualizability results} We will define in Definition \ref{defAF} the 1-morphism $\AA_\eta: \Vect_\FK \to \VV$ induced by the unit inclusion in a braided tensor category $\VV$. It is $\VV$ as a monoidal category with obvious actions. 
Let us recall previously known results about its dualizability. The following is \cite[Theorem 5.1]{GS}, \cite[Theorem 5.16]{BJS} and \cite[Theorem 5.21]{BJS} respectively.
\begin{myth}
The 1-morphism $\AA_\eta$ is always 1-dualizable. It is 2-dualizable as soon as $\VV$ is cp-rigid, and 3-dualizable as soon as $\VV$ is fusion. 
\end{myth}
Note that the requirement fusion can easily be relaxed to rigid finite semisimple, without the assumption that the unit is simple, see the proof of Theorem \ref{thmFTplus1}.

\subsection{Bimodules induced by functors} \label{sectFunBimod}
We give basic definitions and facts about bimodules induced by (braided) monoidal functors, and show how to compute their adjoints.
\subsubsection{Definition and coherence} We show that the notion of bimodules induced by functors behaves as expected in $\BrTens$. Namely, the Morita category, whose morphisms are bimodules, extends the category whose morphisms are functors.
\begin{mydef} \label{defAF} Let $\AA$ and $\BB$ be two objects of $\BrTens$. A braided monoidal functor $F:\AA \to \BB$ induces an $\AA$-$\BB$-central algebra $\AA_F$ which is given by $\BB$ as a monoidal category on which $\AA$ acts on the top using $F(-)\otimes-$ and $\BB$ acts on the bottom using $-\otimes-$. More formally its structure of $\AA$-$\BB$-central algebra is given by:
$$\begin{array}{cll} \AA \boxtimes \BB^{\sigma op} &\to& Z(\BB) \\ (A,B) &\mapsto& (F(A) \otimes B, (\Id_{F(A)} \otimes \sigma^{-1}_{B,-}) \circ (\sigma_{-,F(A)} \otimes \Id_B))
\end{array}$$ where $\sigma$ is the braiding in $\BB$. It is braided monoidal because $F$ is braided monoidal.

It also induces a $\BB$-$\AA$-central algebra $\overline{\AA}_F$ which is also given by $\BB$ as a monoidal category on which $\AA$ acts on the bottom using $-\otimes F(-)$ and $\BB$ acts on the top using $-\otimes-$.

When the functor $F$ is understood, we may write ${}_\AA\BB_\BB$ for $\AA_F$ and ${}_\BB\BB_\AA$ for $\overline{\AA}_F$ 
\end{mydef}
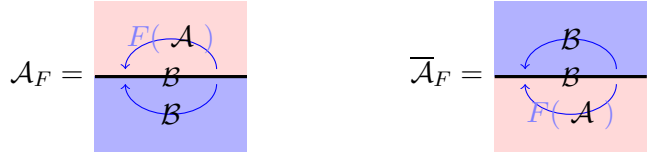
\begin{figure}[h!]
    \centering $\AA_F=$
    \begin{tikzpicture}[scale = 2, baseline = -2pt] 
    \fill[red!15] (0,0) rectangle (1,0.5);
    \fill[blue!30] (0,0) rectangle (1,-0.5);
    \draw[<-, blue] (0.2,0.05) arc (180:0:0.3 and 0.2);
    \draw[<-, blue] (0.2,-0.05) arc (-180:0:0.3 and 0.2);
    \node at (0.5,-0.25){$\BB$};
    \node at (0.5,0.25){\textcolor{blue!50}{$F($} $\AA$ \textcolor{blue!50}{$)$}};
    \draw[very thick] (0,0) -- (1,0) node[midway]{\small $\BB$};
    \end{tikzpicture}\hspace{2cm} $\overline{\AA}_F=$
    \begin{tikzpicture}[scale = 2, baseline = -2pt] 
    \fill[red!15] (0,0) rectangle (1,-0.5);
    \fill[blue!30] (0,0) rectangle (1,0.5);
    \draw[<-, blue] (0.2,0.05) arc (180:0:0.3 and 0.2);
    \draw[<-, blue] (0.2,-0.05) arc (-180:0:0.3 and 0.2);
    \node at (0.5,0.25){$\BB$};
    \node at (0.5,-0.25){\textcolor{blue!50}{$F($} $\AA$ \textcolor{blue!50}{$)$}};
    \draw[very thick] (0,0) -- (1,0) node[midway]{\small $\BB$};
    \end{tikzpicture}
    \caption{The 1-morphisms $\AA_F$ and $\overline{\AA}_F$}
\end{figure}
\begin{prop}The above induced-central-algebra construction preserves composition. Given two braided monoidal functors $F: \AA \to \BB$ and $G:\BB \to \CC$, one has $\AA_G \circ \AA_F \simeq \AA_{G \circ F}$ and $\overline{\AA}_F \circ \overline{\AA}_G \simeq \overline{\AA}_{G\circ F}$.
\end{prop}
\begin{myproof}
We want to prove that ${}_\AA\BB_\BB\underset{\BB}{\boxtimes} {}_\BB\CC_\CC \simeq {}_\AA\CC_\CC$. This is true on the underlying categories as $\BB\underset{\BB}{\boxtimes}\CC \overset{\Phi}{\simeq} \CC$ with equivalence given on pure tensors by $\Phi(B\boxtimes C)= G(B) \otimes C$. This assignment is balanced as $G$ is monoidal: $$\Phi((B\otimes B')\boxtimes C) = G(B\otimes B') \otimes C \simeq G(B) \otimes G(B') \otimes C = \Phi(B\boxtimes (G(B') \otimes C)).$$ 
It is monoidal (the monoidal structure on the relative tensor product is described in \cite[Definition-Proposition 3.6]{BJS}) by:
$$ \Phi(B \boxtimes C) \otimes \Phi(B' \boxtimes C') = G(B) \otimes C \otimes G(B') \otimes C' \underset{\sim}{\overset{\sigma_{C,B'}}{\longrightarrow}} G(B) \otimes G(B') \otimes C \otimes C' \simeq  \Phi((B \boxtimes C) \otimes (B' \boxtimes C')).$$ 
The bottom action of $\CC$ is unchanged, and the top action of $\AA$ is preserved by $\Phi$: 
$$A \rhd (B\boxtimes C) := (A\rhd \unit) \otimes(B\boxtimes C) = (F(A)\otimes B) \boxtimes C \overset{\Phi}{\mapsto} G(F(A)) \otimes G(B) \otimes C = A \rhd \Phi(B\boxtimes C).$$
Finally, let us show that $\Phi$ preserves the central structure.
The central structure in the composed bimodule $\AA_F \underset{\BB}{\boxtimes} \AA_G$ is given by:
$$(B \boxtimes C) \lhd A :=(B \boxtimes C) \otimes (F(A) \boxtimes 1_\CC) = (B \otimes F(A)) \boxtimes C  \underset{\sim}{\overset{\sigma^\BB_{B,F(A)}\boxtimes \Id_C}{\longrightarrow}}(F(A) \otimes B) \boxtimes C = A \rhd (B \boxtimes C)$$
This maps under $\Phi$, using that $G$ is braided monoidal, to $\sigma^\CC_{G(B),G(F(A))}\otimes \Id_C$.
And indeed, the following diagram, where the horizontal arrows are the central structures and the vertical arrow monoidality of $\Phi$, commutes:
\begin{center}
    \begin{tikzpicture}[xscale = 7, yscale=1.5]
\node (D) at (1,0) {$\Phi((A \rhd \unit) \otimes (B \boxtimes C))$};
\node (C) at (0,0) {$\Phi((B \boxtimes C)\otimes( \unit \lhd A))$};
\node (B) at (1,1) {$\Phi(A \rhd \unit) \otimes \Phi(B \boxtimes C)$};
\node (A) at (0,1) {$\Phi(B \boxtimes C)\otimes \Phi( \unit \lhd A)$};
\draw[->] (A)--(B) node[midway, above]{$\sigma^\CC_{(G(B)\otimes C),G(F(A))}$};
\draw[->] (A)--(C) node[midway, left]{$\Id_{G(B)} \otimes \sigma^\CC_{C,G(F(A))}$}; 
\draw[->] (C)--(D) node[midway, above]{$\sigma^\CC_{G(B),G(F(A))}\otimes \Id_C$};
\draw[->] (B)--(D)node[midway, right]{$\Id$};
\end{tikzpicture}
\end{center}
 The $\overline{\AA}$ case is similar.
\end{myproof}
\begin{mydef}\label{defBimodIndFun}
Let $\CC$ and $\DD$ be $\AA$-$\BB$-central algebras, i.e. 1-morphisms of $\BrTens$. A bimodule monoidal functor $F:\CC \to \DD$ preserving the $\AA$-$\BB$-central structures induces an $\AA$-$\BB$-centered $\CC$-$\DD$-bimodule $\MM_F$ which is given by $\DD$ as a category on which $\CC$ acts on the left using $F(-)\otimes-$ and $\DD$ act on the right using $-\otimes-$. The $\AA$-$\BB$-centered structure on $\MM_F$ is induced by the $\AA$-$\BB$-central structure of $\DD$, and the fact that $F$ is a bimodule functor:
$$F(A \rhd \unit_\CC \lhd B) \otimes M \simeq (A \rhd 1_\DD \lhd B) \otimes M \overset{\sigma^\DD}{\simeq} M \otimes (A \rhd 1_\DD \lhd B)$$
It also induces an $\AA$-$\BB$-centered $\DD$-$\CC$-bimodule $\overline{\MM}_F$ which is again given by $\DD$ as a monoidal category on which $\CC$ acts on the right using $-\otimes F(-)$ and $\DD$ act on the left using $-\otimes-$.

When the functor $F$ is understood, we may write ${}_\CC\DD_\DD$ for $\MM_F$ and ${}_\DD\DD_\CC$ for $\overline{\MM}_F$ 
\end{mydef}
\begin{prop}\label{propBimodCompo}
The above induced-bimodule construction preserves:
\begin{enumerate}
    \item horizontal composition:\\ Given two $\AA$-$\BB$-bimodule monoidal functors $F: \CC \to \DD$ and $G:\DD \to \EE$ preserving central structures, one has $\MM_G \circ \MM_F \simeq \MM_{G \circ F}$ and $\overline{\MM}_F \circ \overline{\MM}_G \simeq \overline{\MM}_{G\circ F}$,
    \item vertical composition:\\ Given $\CC$ and $\DD$ two $\AA_1$-$\AA_2$-central algebras, $\CC'$ and $\DD'$ two $\AA_2$-$\AA_3$ central algebras, $F: \CC \to \DD$ an $\AA_1$-$\AA_2$-bimodule monoidal functor and $F':\CC'\to\DD'$ an $\AA_2$-$\AA_3$-bimodule monoidal functor preserving central structures, one has $\MM_F \underset{\AA_2}{\boxtimes} \MM_{F'} \simeq \MM_{F \underset{\AA_2}{\boxtimes} F'}$ and $\overline{\MM}_F \underset{\AA_2}{\boxtimes} \overline{\MM}_{F'} \simeq \overline{\MM}_{F \underset{\AA_2}{\boxtimes} F'}$.
\end{enumerate}
\end{prop}
\begin{myproof}
The first point is similar to the last proposition. We proved that ${}_\CC\DD_\DD\underset{\DD}{\boxtimes} {}_\DD\EE_\EE \overset{\Phi}{\simeq} {}_\CC\EE_\EE$, as bimodules. Recall from \cite[Definition-Proposition 3.13]{BJS} that the centered structure on the composition of bimodules $\DD\underset{\DD}{\boxtimes} \EE$ is given by the composition of the centered structure and a balancing. In our case on some $A,\ D,\ E$, this is: $$D\boxtimes (E \otimes A)\underset{\sim}{\overset{\Id_D \boxtimes\sigma^\EE_{E,A}}{\longrightarrow}} D \boxtimes (A \otimes E) \simeq (D \otimes A) \boxtimes E\underset{\sim}{\overset{\sigma^\DD_{D,A}\boxtimes \Id_E}{\longrightarrow}} (A \otimes D) \boxtimes E $$ which maps by $\Phi$ to $(G(\sigma^\DD_{D,A})\otimes \Id_E)\circ(\Id_{G(D)} \otimes\sigma^\EE_{E,A})$. The centered structure of ${}_\CC\EE_\EE$ is given by $\sigma^\EE_{G(D) \otimes E,A}$. They coincide as $G$ preserves central structures. 

The second point is not surprising either. We want ${}_\CC\DD_\DD \underset{\AA_2}{\boxtimes} {}_{\CC'} {\DD'}_{\DD'} \simeq {}_{\CC\underset{\AA_2}{\boxtimes}\CC'}\DD\underset{\AA_2}{\boxtimes}\DD'_{\DD\underset{\AA_2}{\boxtimes}\DD'}$, which is true on the underlying categories. Because $F$ and $F'$ are bimodule functors, the functor $F \boxtimes F': \CC\boxtimes \CC' \to \DD\boxtimes \DD' \surj \DD\underset{\AA_2}{\boxtimes}\DD'$ is $\BB$-balanced and descends to the relative tensor product $\CC\underset{\AA_2}{\boxtimes}\CC'$. We then see that the left $\CC\underset{\AA_2}{\boxtimes}\CC'$-action is the one induced by $F\boxtimes F'$ on the relative tensor product, namely action by $F\underset{\AA_2}{\boxtimes}F'$. The centered structures are both given by the central structure of $\DD\underset{\AA_2}{\boxtimes}\DD'$ and coincide.
\end{myproof}
\subsubsection{Dualizability}
Given a braided monoidal functor $F:\AA \to \BB$, we will prove that both adjoints of $\AA_F$ are given by $\overline{\AA}_F$. For the right adjunction, the counit should go: $$\AA_F \circ \overline{\AA}_F = {}_\BB\BB_\AA\underset{\AA}{\boxtimes}{}_\AA\BB_\BB \to \Id_\BB = {}_\BB\BB_\BB.$$ We actually have a functor going this way, the tensor product $T$ in $\BB$, which is $\AA$-balanced and descends to the relative tensor product. We denote it $T_{bal}:\BB\underset{\AA}{\boxtimes}\BB \to \BB$, and it is indeed a $\BB$-$\BB$-bimodule monoidal functor. The central structures on both sides are given by braiding in $\BB$, which is preserved by $T$. Hence we can construct a $\BB$-$\BB$-centered ${}_\BB\BB_\AA\underset{\AA}{\boxtimes}{}_\AA\BB_\BB$-${}_\BB\BB_\BB$-bimodule $\MM_{T_{bal}}$ using Definition \ref{defBimodIndFun}.

The unit should go: $$\Id_\AA = {}_\AA\AA_\AA \to \overline{\AA}_F\circ \AA_F = {}_\AA\BB_\BB \underset{\BB}{\boxtimes} {}_\BB \BB_\AA \simeq {}_\AA\BB_\AA.$$ Again we have a functor $F:\AA \to \BB$ which is an $\AA$-$\AA$-module monoidal functor. The central structure on the left is given by braiding in $\AA$, and on the right by braiding in $\BB$. The first is sent on the latter because $F$ is braided monoidal, and the central structures are preserved. Therefore we also have an $\AA$-$\AA$-centered ${}_\AA\AA_\AA$-${}_\AA\BB_\AA$-bimodule $\MM_F$. 

Note also that the identity of $\AA_F$ is the bimodule induced by $\Id_\BB$ seen as an $\AA$-$\BB$-bimodule monoidal functor.

\begin{prop} \label{propAdjAF}
The 1-morphism $\AA_F$ has right adjoint given by $\overline{\AA}_F$, with counit $\MM_{T_{bal}}$ and unit $\MM_F$. Its left adjoint is also given by $\overline{\AA}_F$, with counit $\overline{\MM}_F$ and unit $\overline{\MM}_{T_{bal}}$.
\end{prop}
\begin{myproof}
We directly check the snake. We repeatedly use Proposition \ref{propBimodCompo}:\\
\begin{tikzpicture}[xscale = 3, baseline = 1.3cm]
\node (sIdx) at (0,3) {$\AA$};
\node (tIdx) at (0,1.5) {$\AA$};
\draw[very thin, gray,double] (sIdx)--(tIdx) node[black, midway]{$ {}_\AA\AA_\AA$};
\node (tf1) at (0,0) {$\BB$};
\draw[very thin, gray, ->] (tIdx) -- (tf1) node[black, midway]{$ {}_\AA\BB_\BB$};
\node (seta) at (1,3) {$\AA$};
\node (meta) at (1,2) {$\BB$};
\node (teta) at (1,1) {$\AA$};
\draw[very thin, gray,->] (seta) -- (meta) node[black, midway]{$ {}_\AA\BB_\BB$};
\draw[very thin, gray,->] (meta) -- (teta) node[black, midway]{$ {}_\BB\BB_\AA$};
\node (tf2) at (1,0) {$\BB$};
\draw [very thin, gray,->] (teta) -- (tf2) node[black, midway]{$ {}_\AA\BB_\BB$};
\node (sf) at (2,3){$\AA$};
\node (sIdy) at (2,1.5){$\BB$};
\draw[very thin, gray,->] (sf) -- (sIdy) node[black, midway]{$ {}_\AA\BB_\BB$};
\node (tIdy) at (2,0){$\BB$};
\draw[very thin, gray,double] (sIdy)--(tIdy) node[black, midway]{$ {}_\BB\BB_\BB$};
\draw[very thin, gray!50] (sIdx) -- (seta);
\draw[very thin, gray!50] (tIdx) -- (teta);
\draw[very thin, gray!50] (tf1) -- (tf2);
\draw[very thin, gray!50] (seta) -- (sf);
\draw[very thin, gray!50] (meta) -- (sIdy);
\draw[very thin, gray!50] (tf2) -- (tIdy);
\draw[double, ->] (0.2, 2.2) -- (0.75,2) node[midway, above]{$\MM_F$};
\draw[black, line width = 3pt] (1.3, 2.3) -- (1.7,2.2)  node[midway, above]{$\MM_{\Id_{\BB}}$};
\draw[white, line width = 2pt] (1.25, 2.31) -- (1.75,2.19);
\draw[double, ->] (1.2, 0.8) -- (1.8,0.6) node[midway, above]{$\MM_{T_{bal}}$};
\draw[black, line width = 3pt] (0.3, 0.4) -- (0.7,0.3) node[midway, above]{$\MM_{\Id_{\BB}}$};
\draw[white, line width = 2pt] (0.25, 0.41) -- (0.75,0.29);
\node (Fsf) at (-1,2.3) {$\AA$};
\node (Ftf) at (-1,0.7) {$\BB$};
\draw [very thin, gray,->] (Fsf) -- (Ftf) node[black, midway]{$ {}_\AA\BB_\BB$};
\node (Lsf) at (3,2.3) {$\AA$};
\node (Ltf) at (3,0.7) {$\BB$};
\draw [very thin, gray,->] (Lsf) -- (Ltf) node[black, midway]{$ {}_\AA\BB_\BB$};
\draw[very thin, gray!50] (Fsf) -- (sIdx);
\draw[very thin, gray!50] (Ftf) -- (tf1);
\draw[very thin, gray!50] (sf) -- (Lsf);
\draw[very thin, gray!50] (tIdy) -- (Ltf);
\node at (-0.5,1.5){$\simeq$};
\node at (2.5,1.5){$\simeq$};
\node (3') at (0.5,-0.7) {\LARGE $\simeq$};
\end{tikzpicture}
\\
\begin{tikzpicture}[xscale = 3, baseline = 1.3cm]
\node (sIdx) at (0,3) {$\AA$};
\node (tIdx) at (0,1.5) {$\AA$};
\draw[very thin, gray,double] (sIdx)--(tIdx) node[black, midway]{$ {}_\AA\AA_\AA$};
\node (tf1) at (0,0) {$\BB$};
\draw[very thin, gray, ->] (tIdx) -- (tf1) node[black, midway]{$ {}_\AA\BB_\BB$};
\node (seta) at (1,3) {$\AA$};
\node (meta) at (1,2) {$\BB$};
\node (teta) at (1,1) {$\AA$};
\draw[very thin, gray,->] (seta) -- (meta) node[black, midway]{$ {}_\AA\BB_\BB$};
\draw[very thin, gray,->] (meta) -- (teta) node[black, midway]{$ {}_\BB\BB_\AA$};
\node (tf2) at (1,0) {$\BB$};
\draw [very thin, gray,->] (teta) -- (tf2) node[black, midway]{$ {}_\AA\BB_\BB$};
\node (sf) at (2,3){$\AA$};
\node (sIdy) at (2,1.5){$\BB$};
\draw[very thin, gray,->] (sf) -- (sIdy) node[black, midway]{$ {}_\AA\BB_\BB$};
\node (tIdy) at (2,0){$\BB$};
\draw[very thin, gray,double] (sIdy)--(tIdy) node[black, midway]{$ {}_\BB\BB_\BB$};
\draw[very thin, gray!50] (sIdx) -- (seta);
\draw[very thin, gray!50] (tf1) -- (tf2);
\draw[very thin, gray!50] (seta) -- (sf);
\draw[very thin, gray!50] (tf2) -- (tIdy);
\draw[double, ->] (0.25, 1) -- (0.75,1) node[midway, above]{$\MM_{F\underset{\AA}{\boxtimes}\Id_\BB}$};
\draw[double, ->] (1.25, 1) -- (1.75,1) node[midway, above]{$\MM_{\Id_\BB\underset{\BB}{\boxtimes}T_{bal}}$};
\node (Fsf) at (-1,2.3) {$\AA$};
\node (Ftf) at (-1,0.7) {$\BB$};
\draw [very thin, gray,->] (Fsf) -- (Ftf) node[black, midway]{$ {}_\AA\BB_\BB$};
\node (Lsf) at (3,2.3) {$\AA$};
\node (Ltf) at (3,0.7) {$\BB$};
\draw [very thin, gray,->] (Lsf) -- (Ltf) node[black, midway]{$ {}_\AA\BB_\BB$};
\draw[very thin, gray!50] (Fsf) -- (sIdx);
\draw[very thin, gray!50] (Ftf) -- (tf1);
\draw[very thin, gray!50] (sf) -- (Lsf);
\draw[very thin, gray!50] (tIdy) -- (Ltf);
\node at (-0.5,1.5){$\MM_{\unit_\AA\boxtimes-}$};
\node at (2.5,1.5){$\MM_{-\otimes-}$};
\end{tikzpicture}
\\
which is the bimodule induced by the composition:\\
\begin{tikzpicture}[xscale = 3]
\node (1) at (-1,-1) {$\BB$};
\node (2) at (0,-1) {$\AA\underset{\AA}{\boxtimes}\BB$};
\node (3) at (1,-1) {$\BB\underset{\BB}{\boxtimes}\BB\underset{\AA}{\boxtimes}\BB$};
\node (4) at (2,-1) {$\BB\underset{\BB}{\boxtimes}\BB$};
\node (5) at (3,-1) {$\BB$};
\draw[->] (1) -- (2); \draw[->] (2) -- (3); \draw[->] (3) -- (4); \draw[->] (4) -- (5);
\node (1') at (-1,-2) {$X$};
\node (2') at (0,-2) {$(\unit_\AA,X)$};
\node (3') at (1,-2) {$(\unit_\BB, \unit_\BB,X)$};
\node (4') at (2,-2) {$(\unit_\BB,X)$};
\node (5') at (3,-2) {$X$};
\node at (-0.5,-2){$\mapsto$};\node at (0.5,-2){$\mapsto$};\node at (1.5,-2){$\mapsto$};\node at (2.5,-2){$\mapsto$};
\end{tikzpicture}\\
which is indeed the identity.

Every other snake identity is very similar, with functors going in the other direction for the left adjunction.
\end{myproof}
\begin{prop} \label{propRAforMF}
Let $F:\CC\to\DD$ be an $\AA$-$\BB$-bimodule monoidal functor. The bimodule $\MM_F$ has right adjoint given by $\overline{\MM}_F$, with counit $T_{bal}: \DD\underset{\CC}{\boxtimes}\DD \to \DD$ seen as a $\DD$-$\DD$-bimodule functor and unit $F$ seen as a $\CC$-$\CC$-bimodule functor.
\end{prop}
\begin{myproof}
The proof is the same as above, except that the horizontal morphisms are now the functors instead of the bimodules induced by the functors. The snake identities read:
\begin{align}\label{equSnakeMF}
    (\Id_\DD\underset{\DD}{\boxtimes}T_{bal})\circ (F\underset{\CC}{\boxtimes}\Id_\DD)  \simeq \Id_{\MM_F} \quad \text{and}\quad (T_{bal}\underset{\DD}{\boxtimes} \Id_\DD)\circ (\Id_\DD \underset{\CC}{\boxtimes}F)  \simeq \Id_{\overline{\MM}_F}
\end{align}  as has been used above. Here $\Id_\DD$ is seen alternatively as a $\CC$-$\DD$-bimodule functor and as a $\DD$-$\CC$-bimodule functor.
\end{myproof}
We would like to apply Proposition \ref{propRedundancyDD}.1, to have the left adjoint of $\MM_F$. We need $F$ and $T_{bal}$ to have right adjoints in $\BrTens$. There is a well-known sufficient condition for this.
\begin{prop}[Proposition 4.2 and Corollary 4.3 in {\cite{BJS}}] \label{propRA3morph}
Let $F: \MM \to \NN$ be an $\AA$-$\BB$-centered $\CC$-$\DD$-bimodule functor, so a 3-morphism in $\BrTens$. Suppose that $\MM$ and $\NN$ have enough cp, that $\AA,\ \BB,\ \CC,\ \DD$ are cp-rigid, and that $F$ preserves cp. Then $F^R:\NN \to\MM$ is an $\AA$-$\BB$-centered $\CC$-$\DD$-bimodule functor, and is the right adjoint of $F$ in $\BrTens$.
\end{prop}
All we need to check is that both $F$ and $T_{bal}$ preserve cp.
\begin{mylemma}\label{lemmacpRelTens}
Let $\MM$ and $\NN$ be right and left modules over $\CC$ and $F:\MM\boxtimes \NN \to \PP$ be a cocontinuous $\CC$-balanced functor. Suppose $\MM$ and $\NN$ have enough cp, $\CC$ is cp-rigid and $F$ preserves cp. Then the induced functor $F_{bal}: \MM\underset{\CC}{\boxtimes}{\NN}\to\PP$ preserves cp.

In particular, if $\AA$ and $\BB$ are cp-rigid, then $T_{bal}: \BB\underset{\AA}{\boxtimes}\BB \to \BB$ preserves cp.
\end{mylemma}
\begin{myproof}
Following the proof of closure under composition of 1-morphisms \cite[Section 4.2]{BJS}, the cp objects of $\MM\underset{\CC}{\boxtimes}{\NN}$ are generated by pure tensors of cp objects. These are sent to cp objects in $\PP$.

For the second point, $T_{bal}$ is induced by $T$ which preserves cp as $\BB$ is cp-rigid.
\end{myproof}
We can summarize the result as follows:
\begin{prop}\label{propLAforMF}
Let $F:\CC\to\DD$ be an $\AA$-$\BB$-bimodule monoidal functor which preserves cp, where $\AA,\ \BB,\ \CC,\ \DD$ are cp-rigid. The bimodule $\MM_F$ has left adjoint given by $\overline{\MM}_F$, with counit $F^R$ seen as a $\CC$-$\CC$-bimodule functor and unit ${T_{bal}^R}$ seen as a $\DD$-$\DD$-bimodule functor.
\end{prop}

\subsection{Unit inclusion} We give explicitly the dualizability data of the 1-morphism induced by the unit inclusion in a braided tensor category $\VV$, and criteria for dualizability when $\VV$ has enough cp.
\begin{mydef}
Let $\VV\in \BrTens$ be an $\mathbb E_2$-algebra in $\Pr$. We denote by $T:\VV\boxtimes\VV \to \VV$ its monoidal structure, and $\eta: \Vect_\FK \to\VV$ the inclusion of the unit. The functor $\eta$ is braided monoidal and induces a $\Vect_\FK$-$\VV$-central algebra $\AA_\eta$, namely a 1-morphism in $\BrTens$. Remember that we denote by $\Aflat_\eta \in \BrTens^\to$ the associated object in the oplax arrow category.
\end{mydef}
\begin{myth}\label{thmVcpgen}
The 1-morphism $\AA_\eta$ is both twice left and twice right adjunctible, with adjunctibility data as displayed:
\begin{center} \begin{tikzpicture}[baseline = 10pt, xscale = 1.4, yscale = 1.2]
\node (A) at (2,2) {$\AA_\eta$};
\node[inner sep = 20pt] (L) at (0,1) {$\vert$};
\node[inner sep = 20pt] (R) at (4,1) {$\vert$};
\node[left of= R] (Rcounit)  {$\MM_T$};
\node[right of= R] (Runit){$\MM_\eta$};
\node[left of= L] (Lcounit) {$\overline{\MM}_\eta$};
\node[right of= L] (Lunit) {$\overline{\MM}_T$};
\node (LLc) at (-2,0) {$T\vert \eta$};
\node (LLu) at (0.5,0) {$T_{bal}\vert T$};
\node (RRc) at (3.5,0) {$T_{bal}\vert T$};
\node (RRu) at (6,0) {$T\vert \eta$};
\draw[->] (A) -- (R) node[near end,above right]{\footnotesize $\overline{\AA}_\eta$};
\draw[->] (A) -- (L) node[near end,above left]{\footnotesize $\overline{\AA}_\eta$};
\draw[->] (Rcounit) -- (RRc) node[near end,above right]{\tiny $\overline{\MM}_T$};
\draw[->] (Runit) -- (RRu) node[near end,above right]{\tiny $\overline{\MM}_\eta$};
\draw[->] (Lcounit) -- (LLc) node[near end,above left]{\tiny $\MM_\eta$};
\draw[->] (Lunit) -- (LLu) node[near end,above left]{\tiny $\MM_T$};
\end{tikzpicture}\end{center}
where $T_{bal}: \VV \underset{\VV \boxtimes \VV}{\boxtimes} \VV \to \VV$ is induced by $T$ on the relative tensor product
\end{myth}
\begin{myproof}
We use the results of Section \ref{sectFunBimod}. By Proposition \ref{propAdjAF}, the 1-morphism $\AA_\eta$ has left and right adjoints given by $\overline{\AA}_\eta$, with units and counits as displayed in the second line above, with $\eta: \Vect_\FK \to \VV$ now seen as a $\Vect_\FK$-$\Vect_\FK$-bimodule monoidal functor, and $T: \VV \underset{\Vect_\FK}{\boxtimes}\VV\to\VV$ the tensor product balanced over $\Vect_\FK$ so not balanced.  

Then by Proposition \ref{propRAforMF} each of these bimodules has either a left or a right adjoint, with units and counits as displayed, with $T_{bal}: \overline{\MM}_T \underset{\VV\boxtimes\VV}{\boxtimes} \MM_T = \VV\underset{\VV\boxtimes\VV}{\boxtimes}\VV \to \VV$ induced by $T$.
\end{myproof}

\begin{corr}\label{corrUnitDDoplax}
The object $\Aflat_\eta$ is 2-dualizable in $\BrTens^\to$, and:\\
$Ru(Ru(\Aflat_\eta))$ has a right adjoint if and only if both $T_{bal}$ and $Ru(Ru(\VV))$ do. \\
$Rco(Ru(\Aflat_\eta))$ has a right adjoint if and only if both $T$ and $Rco(Ru(\VV))$ do. \\
$Rco(Rco(\Aflat_\eta))$ has a right adjoint if and only if both $\eta$ and $Rco(Rco(\VV))$ do.
\end{corr}
\begin{myproof} For 2-dualizability, we use the criterion of \cite[Theorem 7.6]{JFS}, we know that $\VV$ is 2-dualizable by \cite[Theorem 5.1]{GS} and $\AA_\eta$ is twice right adjunctible by the theorem above. The rest is Theorem \ref{thmDDinOplax} on the right dualizability data of $\AA_\eta$.
\end{myproof}
\begin{myth}\label{thmCrit3oplax} Suppose that $\VV$ has enough cp, then $\Aflat_\eta$ is 3-dualizable if and only if $\VV$ the free cocompletion of a small rigid braided monoidal category.
\end{myth}
\begin{myproof} The heart of the proof is to notice that $T$ appears in the dualizability data, and by \cite[Proposition 4.1]{BJS} when $\VV$ has enough cp, it is cp-rigid if and only if $T$ has a bimodule cocontinuous right adjoint.

If $\Aflat_\eta$ is 3-dualizable then $Ru(Ru(\Aflat_\eta))$, $Rco(Ru(\Aflat_\eta))$ and $Rco(Rco(\Aflat_\eta))$ have right adjoints, so $T_{bal}$, $T$ and $\eta$ have bimodule cocontinuous right adjoints. The functors $T$ and $\eta$ preserving cp mean that they are well-defined on $\mathscr{V}:=\VV^{cp}$ and endow it with a monoidal structure, and $\Vv$ is rigid as $\VV$ is cp-rigid. Therefore $\VV$ is the free cocompletion of a small rigid braided monoidal category.

On the other hand if $\VV$ is the free cocompletion of a small rigid braided monoidal category then it is cp-rigid and hence 3-dualizable, \cite[Theorem 5.16]{BJS}. The functors $T$ and $\eta$, and also $T_{bal}$ by Lemma \ref{lemmacpRelTens}, preserve cp, and have bimodule cocontinuous right adjoints by Proposition \ref{propRA3morph}. We get that $\AA_\eta$ is 3-times right adjunctible and its source and targets are 3-dualizable, so $\Aflat_\eta$ 3-dualizable by \cite[Theorem 7.6]{JFS}. 
\end{myproof}
\begin{myth}\label{thmCritnc3d} Suppose that $\VV$ has enough cp, then $\Aflat_\eta$ is non-compact-3-dualizable if and only if $\VV$ is cp-rigid.
\end{myth}
\begin{myproof}
If $\VV$ is cp-rigid, then $\VV$ is 3-dualizable and $T$ and $T_{bal}$ have right adjoints in $\BrTens$. By Corollary \ref{corrUnitDDoplax}, $\Aflat_\eta$ is non-compact-3-dualizable.

Suppose now $\Aflat_\eta$ non-compact-3-dualizable, then $T$ has a bimodule cocontinuous right adjoint, and $\VV$ is cp-rigid.
\end{myproof}
\begin{myth}\label{thmFTplus1}
Let $\VV$ be a braided tensor category with enough cp. Then the following are equivalent: \begin{enumerate}
    \item $\AA_\eta$ is 3-dualizable,
    \item $\AA_\eta$ is 3-adjunctible, and
    \item $\VV$ is rigid finite semisimple.
\end{enumerate}
\end{myth}
\begin{myproof}
The implication $1\Rightarrow 2$ is immediate: for a 1-morphism 3-dualizable demands 3-adjunctible and 4-dualizablility of the source and target.

The implication $3\Rightarrow 1$ is essentially  \cite[Theorem 5.21]{BJS}. If $\VV$ is fusion, then $\VV$ and $\AA_\eta$ lie $BrFus$ which has duals. Now fusion demands simplicity of the unit, which may not be the case here. This is easily solved by noticing that coproduct agrees with product in $\Pr$ and ought to be called direct sum \cite[Remark 2.5]{BCJReflDualPr}, and that braided rigid finite semisimple categories are direct sums of fusion categories \cite[Section 4.3]{EGNO}.

Let us prove $2\Rightarrow 3$. If $\AA_\eta$ is 3-adjunctible then $\MM_\eta$ and $\overline{\MM}_\eta$, which are respectively $Ru(\AA_\eta)$ and $Lco(\AA_\eta)$ by Theorem \ref{thmCrit3oplax}, must be 2-adjunctible. Hence their composite $\MM_\eta \underset{\VV}{\boxtimes}\overline{\MM}_\eta $ has to be 2-adjunctible in the symmetric monoidal 2-category $\Omega\Omega\BrTens\simeq \Pr$. This composition is just $\VV\underset{\VV}{\boxtimes}\VV\simeq\VV$ as a category, and by our assumption that it has enough cp, it actually lies in the full subcategory $\Bimod_\FK\subseteq \Pr$. By \cite[Theorem A.22]{BDSPV}, the 2-dualizable objects of $\Bimod_\FK$ are finite semisimple categories. We already saw that $\VV$ has to be cp-rigid, so $\VV^{cp}$ is rigid finite semisimple, and so is $\VV\simeq \Free(\VV^{cp})$.
\end{myproof}
\begin{rmk}
A very similar result one categorical dimension down, in $\Alg_1(Rex_\Cc)$, is proven in \cite[Theorem B]{FreedTelemanGapped}. The proof is similar too, but we couldn't directly use their result on $\MM_\eta$ as we work in $\Bimod_\FK$ and not in $Rex_\Cc$.
\end{rmk}
\begin{rmk}\label{rmkNon4DV}
Both results need full adjunctibility of $\AA_\eta$: oplax dualizability does not imply semisimplicity, one can take the free cocompletion of a non-semisimple ribbon category in Theorem \ref{thmCrit3oplax}. Semisimplicity is not needed for 4-dualizability either, as proven in \cite{BJSS}. However, if we assume that $\VV$ is 4-dualizable and $\AA_\eta$ is 3-oplax-dualizable, which is the case of interest for Section \ref{sectAppli}, then to-appear work of William Stewart shows that $\AA_\eta$ is 3-adjunctible. This has an interesting consequence: the free cocompletion of a ribbon category which is not semisimple cannot be 4-dualizable. Indeed if it were Stewart's result would apply and $\VV$ would have to be semisimple. This justifies that, given a non-semisimple ribbon tensor category as in \cite{CGHP}, we want to work with its Ind-completions, and not its free cocompletion.
\end{rmk}
\begin{rmk}\label{rmkAdjImpliesDual} Being dualizable for a morphism is both a condition on its adjunctibility and on the dualizability of its source and target. However, we saw in the proof of Theorem \ref{thmCrit3oplax} that $\AA_\eta \text{ 3-right-adjunctible }\ \Leftrightarrow\ \AA_\eta$ 3-oplax-dualizable, and in the theorem above that $\AA_\eta \text{ 3-adjunctible }\ \Leftrightarrow\ \AA_\eta$ 3-dualizable. This phenomenon seems to be specific to the unit inclusion.
\end{rmk}
\begin{prop}\label{propVcprigid} 
Suppose that $\VV$ is cp-rigid, then $\AA_\eta$ is 2-adjunctible with the following adjunctibility data in $\BrTens$:
\begin{center} \begin{tikzpicture}[baseline = 10pt, xscale = 1.4, yscale = 1.2]
\node (A) at (2,2) {$\AA_\eta$};
\node[inner sep = 20pt] (L) at (0,1) {$\vert$};
\node[inner sep = 20pt] (R) at (4,1) {$\vert$};
\node[left of= R] (Rcounit)  {$\MM_T$};
\node[right of= R] (Runit){$\MM_\eta$};
\node[left of= L] (Lcounit) {$\overline{\MM}_\eta$};
\node[right of= L] (Lunit) {$\overline{\MM}_T$};
\node (LLc) at (-2,0) {$T\vert \eta$};
\node (RLc) at (-0.5,0) {$\widetilde{\eta^R}\vert T^R$};
\node (LLu) at (0.5,0) {$T_{bal}\vert T$};
\node (RLu) at (1.5,0) {$T^R\vert {T_{bal}^R}$};
\node (LRc) at (2.5,0) {$T^R\vert {T_{bal}^R}$};
\node (RRc) at (3.5,0) {$T_{bal}\vert T$};
\node (LRu) at (4.5,0) {$\widetilde{\eta^R}\vert T^R$};
\node (RRu) at (6,0) {$T\vert \eta$};
\draw[->] (A) -- (R) node[near end,above right]{\footnotesize $\overline{\AA}_\eta$};
\draw[->] (A) -- (L) node[near end,above left]{\footnotesize $\overline{\AA}_\eta$};
\draw[->] (Rcounit) -- (LRc) node[near end,above left]{\tiny $\overline{\MM}_T$};
\draw[->] (Rcounit) -- (RRc) node[near end,above right]{\tiny $\overline{\MM}_T$};
\draw[->] (Runit) -- (LRu) node[near end,above left]{\tiny $\overline{\MM}_\eta$};
\draw[->] (Runit) -- (RRu) node[near end,above right]{\tiny $\overline{\MM}_\eta$};
\draw[->] (Lcounit) -- (LLc) node[near end,above left]{\tiny $\MM_\eta$};
\draw[->] (Lcounit) -- (RLc) node[near end,above right]{\tiny $\MM_\eta$};
\draw[->] (Lunit) -- (LLu) node[near end,above left]{\tiny $\MM_T$};
\draw[->] (Lunit) -- (RLu) node[near end,above right]{\tiny $\MM_T$};
\end{tikzpicture}\end{center}
where $\widetilde{\eta^R}$ is the essentially unique cocontinuous functor that agrees with $\eta^R$ on cp objects.
\end{prop}
\begin{myproof}
The snake for $T^R$ and $\widetilde{\eta^R}$ comes from the following. Denote $\mathscr{V}= \VV^{cp}$. 
Then $T^R$ is computed as the coend $T^R(\unit_\VV) = \displaystyle \int^{(V, W)\in\mathscr{V}^{\otimes 2}} (V \boxtimes W) \otimes \Hom_\VV(V\otimes W,\unit_\VV) \simeq \int^{V\in\mathscr{V}} V \boxtimes V^*$, and more generally $T^R(X)\simeq \displaystyle \int^{V\in\mathscr{V}} (X\otimes V) \boxtimes V^*\simeq \int^{V\in\mathscr{V}} V \boxtimes (V^*\otimes X)$.
For $X$ cp, the snake goes 
$$\begin{array}{cl}
    (\widetilde{\eta^R} \underset{\Vect_\FK}{\boxtimes} \Id_\VV) \circ (\Id_\VV \underset{\VV}{\boxtimes} T^R) (X)& \simeq \displaystyle \int^{V\in\mathscr{V}} \widetilde{\eta^R}(X\otimes V) \boxtimes V^*\\
    &\displaystyle= \int^{V\in\mathscr{V}} \Hom(\unit_\VV,X\otimes V) \otimes V^*\\
    &\displaystyle\simeq \int^{V\in\mathscr{V}} \Hom(V^*,X) \otimes V^* \simeq X
\end{array}$$
The part with $T^R$ and ${T_{bal}^R}$ is given by Proposition \ref{propLAforMF}. Indeed $T$, and hence $T_{bal}$, preserves cp as $\VV$ is cp-rigid.
\\
The fact that this is sufficient for 2-adjunctibility is \cite[Lemma 7.11]{JFS}.
\end{myproof}
\begin{rmk} We studied the oplax-dualizability of $\AA_\eta$ above, but \cite{JFS} also define a notion of lax-dualizability. We are interested in the oplax-dualizability for our applications, but let us include the lax version of our results.
By Theorem \ref{thmVcpgen} $\AA_\eta$ is always 2-lax-dualizable, and it is 3-times left adjunctible if and only if $\eta,\ T$ and $T_{bal}$ have left adjoints in $\BrTens$. 
Using the proposition above, we can also get another characterisation of adjunctibility: every morphism appearing there must have a right adjoint. If $\VV$ has enough cp, then $\AA_\eta$ is 3-adjunctible if and only if $\VV$ is cp-rigid and $\eta,\ \eta^R,\ T^R$ and $T_{bal}^R$ preserve cp.
\end{rmk}
We studied the unit inclusion, but similar arguments work for any bimodule induced by a functor. Instead of a necessary and sufficient condition, we only have a sufficient condition because $T$ no longer appears in the dualizability data, only some balanced version does.
\begin{myth}\label{thm3oplaxAnyBrMonFun}
Let $F: \VV \to \WW$ be a braided monoidal functor between two objects of $\BrTens$. Then the object $\Aflat_F \in \BrTens^\to$ induced by the 1-morphism $\AA_F$ is 2-dualizable. It is non-compact-3-dualizable as soon as $\VV$ and $\WW$ are cp-rigid. In this case, it is 3-dualizable if and only if $F$ preserves cp.
\end{myth}
\begin{myproof}
We know that $Radj(\AA_F) = \overline{\AA}_F$ with $Ru(\AA_F) = \MM_F$ and $Rco(\AA_F)=\MM_{T_{\VV-bal}}$ by Proposition \ref{propAdjAF}, where $T_{\VV-bal}: \WW \underset{\VV}{\boxtimes}\WW \to \WW$ is induced by the monoidal structure on $\WW$.

Then, $Radj(\MM_{T_{\VV-bal}}) = \overline{\MM}_{T_{\VV-bal}}$ with $Ru(\MM_{T_{\VV-bal}}) = {T_{\VV-bal}}$ and $Rco(\MM_{T_{\VV-bal}}) = T_{2bal}$ by Proposition \ref{propRAforMF}, where $T_{2bal}: \WW\underset{\WW \underset{\VV}{\boxtimes}\WW}{\boxtimes}\WW \to \WW$ is induced by the monoidal structure on $\WW$.

Similarly, $Radj(\MM_F) = \overline{\MM}_F$ with $Ru(\MM_F) = F$ and $Rco(\MM_F)=T_{\VV-bal}$.

We know by Theorem \ref{thmDDinOplax} that the existence and right adjunctibility of $Ru(Ru(\Aflat_F))$, $Rco(Ru(\Aflat_F))$ and $Rco(Rco(\Aflat_F))$ is equivalent to that of respectively $T_{2bal}$, $T_{\VV-bal}$ and $F$, and of the same units/counits of the source and target. So $\Aflat_F$ is non-compact-3-dualizable if and only if $T_{\VV-bal}$ and $T_{2bal}$ have right adjoints in $\BrTens$, and both $\VV$ and $\WW$ are non-compact-3-dualizable. This is true as soon as $\VV$ and $\WW$ are cp-rigid by Lemma \ref{lemmacpRelTens} and \cite[Theorem 5.6]{BJS}.

It is 3-dualizable if and only if $F$, $T_{\VV-bal}$ and $T_{2bal}$ have right adjoints and $\VV$ and $\WW$ are 3-dualizable. If $\VV$ and $\WW$ are cp-rigid, this is true if and only if $F$ preserves cp.
\end{myproof}

\subsection{The relative theory on the circle} \label{subsCircle}
We compute the value on the circle of the relative TQFT $\RR_\VV$ induced by $\Aflat_\eta$ under the cobordism hypothesis, for any $\VV$. Namely, we write $S^1_{nb} = ev_{pt}\circ coev_{pt}$, compute the images of $ev_{pt}$ and $coev_{pt}$ under $\RR_\VV$, which are $ev_{\Aflat_\eta}$ and $coev_{\Aflat_\eta}$, and compose them. Note that it is $S^1$ with non-bounding framing that we are computing\footnote{In dimension 3, there are two framings on the circle, only one of which bounds a framed disk.}. We need the symmetric monoidal structure of $\CC$ to compose $ev_X:\unit\to X\otimes X^*$ and $coev_X: X\otimes X^*\simeq X^*\otimes X\to \unit$. We know that the evaluation and coevaluation for $\Aflat_\eta$ are mates of the unit and counit for the right adjunction of $\AA_\eta$, namely $\MM_\eta$ and $\MM_T$. It might sound surprising that one can compose them, but indeed up to whiskering and mating they are composeable, see Figure \ref{unitCounitCompose}.
\begin{figure}[h!]
    \centering
\begin{tikzpicture}[scale = 1.2]
\fill[gray!20, fill opacity = 0.7] (0,-1)arc (-90:90:0.5 and 1) -- (2,1)..controls (1.7,1) and (1.5,0.5)..(1,0).. controls (2,-1.5) and (3,1.5) .. (4,0)..controls (3.5,-0.5) and (3.3,-1)..(3,-1) -- (0,-1);
\fill[gray!20, fill opacity = 0.3] (0,-1)arc (270:90:0.5 and 1) -- (2,1).. controls (2.3,1) and (2.7,-1)..(3,-1) -- cycle;
\fill[blue!40, fill opacity = 0.3] (5,-1) arc (270:90:0.5 and 1) -- (2,1).. controls (2.3,1) and (2.7,-1)..(3,-1) -- cycle;
\fill[blue!40, fill opacity = 0.7] (5,-1) arc (-90:90:0.5 and 1) -- (2,1)..controls (1.7,1) and (1.5,0.5)..(1,0).. controls (2,-1.5) and (3,1.5) .. (4,0)..controls (3.5,-0.5) and (3.3,-1)..(3,-1)--cycle;
\draw (0,-1) arc (270:90:0.5 and 1);
\draw (0,-1) arc (-90:90:0.5 and 1);
\draw[dashed] (5,-1) arc (270:90:0.5 and 1);
\draw (5,-1) arc (-90:90:0.5 and 1);
\draw (0,1) -- (2,1)--(5,1);
\draw (0,-1) -- (3,-1)--(5,-1);
\node at (1,0) {$\bullet$};
\node at (0.8,0.25) {$\MM_\eta$};
\node at (4,0) {$\bullet$};
\node at (4.15,-0.25) {$\MM_T$};
\draw (1,0)..controls (1.5,0.5) and (1.7,1)..(2,1)node[pos = 0.5]{${\AA}_\eta$};
\draw[dashed] (2,1).. controls (2.3,1) and (2.7,-1)..(3,-1) node[black!70, pos = 0.67, sloped, above = -3pt]{\small mating};
\draw (3,-1)..controls (3.3,-1) and (3.5,-0.5).. (4,0);
\draw (1,0).. controls (2,-1.5) and (3,1.5) .. (4,0) node[pos = 0.2]{$\overline{\AA}_\eta$};
\node[gray] at (0,0){$\Vect_\FK$};
\node[scale = 1.2] at (5,0){\large $\VV$};
\end{tikzpicture}    
\caption{The unit and the counit compose up to mating. Beware that the framing is not faithfully represented in this picture.}
    \label{unitCounitCompose}
\end{figure}

We know from \cite[Theorem 2.19]{BJSS} that the evaluation and coevaluation for $\VV$ are respectively ${}_{\VV \boxtimes\VV^{\tiny \sigma op}}\VV_{\Vect_\FK}$ and ${}_{\Vect_\FK}\VV_{\VV^{\tiny \sigma op} \boxtimes\VV}$. Then, Example \ref{exmplDualinCto} gives: $$\RR_\VV(ev_{pt}) =
\begin{tikzpicture}[baseline = 10pt,xscale = 3,yscale = 1.5]
\node (s1) at (0,1) {$\Vect_\FK$};
\node (t1) at (1,1) {$\VV \otimes \VV^{\sigma op}$};
\node (s2) at (0,0) {$\Vect_\FK$};
\node (t2) at (1,0) {$\Vect_\FK$};
\draw[->] (s1) -- (t1) node[midway, above=-2pt]{$\AA_\eta \otimes (\overline{\AA}_\eta)^*$};
\draw[->] (s2) -- (t2) node[midway, below=-2pt]{$\Vect_\FK$};
\draw[->] (s1) -- (s2) node[midway, left=-2pt]{$\Id$};
\draw[->] (t1) -- (t2) node[midway, right=-2pt]{$ev_\VV$};
\draw[double, ->] (0.25,0.25) -- (0.75,0.75)  node[midway, below right=-2pt]{$\MM_\eta$};
\end{tikzpicture}\text{ and } \RR_\VV(coev_{pt}) =
\begin{tikzpicture}[baseline= 10pt, scale = 2]
\node (s1) at (0,1) {$\Vect_\FK$};
\node (t1) at (2.1,1) {$\Vect_\FK$};
\node (s2) at (0,0) {$\Vect_\FK$};
\node (t2) at (2.1,0) {$\VV \otimes \VV^{\sigma op}$};
\draw[->] (s1) -- (t1) node[midway, above=-2pt]{$\Vect_\FK$};
\draw[->] (s2) -- (t2) node[midway, below=-2pt]{$\AA_\eta \otimes (\overline{\AA}_\eta)^*$};
\draw[->] (s1) -- (s2) node[midway, left=-2pt]{$\Id$};
\draw[->] (t1) -- (t2) node[midway, right=-2pt]{$coev_\VV$};
\draw[double, ->] (0.2,0.2) -- (1.75,0.75)  node[sloped, pos=0, below right=-2pt]{\footnotesize $(\MM_T\otimes \Id_{\VV^{\tiny \sigma op}}) \circ_1 \Id_{coev_\VV}$};
\end{tikzpicture}$$
Their composition is vertical stacking and gives that $\RR_\VV(S^1_{nb})$ is the following composition. The blue lines give the connection with Figure \ref{unitCounitCompose}, with correct framing:
\begin{figure}[h!]
    \centering
\begin{tikzpicture}[xscale = 1.3, yscale = 1.2]
\node (V1) at (0,4) {$\Vect_\FK$};
\node (V2) at (1.5,4) {$\Vect_\FK$};
\node (V3) at (4,4) {$\Vect_\FK$};
\node (V4) at (7,4) {$\Vect_\FK$};
\node (V5) at (10,4) {$\Vect_\FK$};
\node (W1) at (0,0) {$\Vect_\FK$};
\node (W2) at (1.5,0) {$\Vect_\FK$};
\node (W3) at (4,0) {$\Vect_\FK$};
\node (W4) at (7,0) {$\Vect_\FK$};
\node (W5) at (10,0) {$\Vect_\FK$};
\draw[double] (V1)--(V2)--(V3)--(V4)--(V5);
\draw[double] (W1)--(W2)--(W3)--(W4)--(W5);
\node (A2) at (1.5,2) {$\VV$};
\node (A3) at (3.5,3) {$\VV$};
\node (A'3) at (4.5,3) {$\Vect_\FK$};
\node (AT3)at (4,3) {$\otimes$};
\node (B3) at (3.2,2) {$\VV$};
\node (B'3) at (4,2) {$\VV^{\tiny \sigma op}$};
\node (B''3) at (4.8,2) {$\VV$};
\node (BT3) at (3.6,2) {$\otimes$};
\node (BT'3) at (4.4,2) {$\otimes$};
\node (C3) at (3.5,1) {$\Vect_\FK$};
\node (C'3) at (4.5,1) {$\VV$};
\node (CT3)at (4,1) {$\otimes$};
\node (A4) at (6.5,3) {$\VV$};
\node (A'4) at (7.5,3) {$\VV^{\tiny \sigma op}$};
\node (AT4) at (7,3) {$\otimes$};
\node (B4) at (6.5,2) {$\Vect_\FK$};
\node (C4) at (6.5,1) {$\VV$};
\node (C'4) at (7.5,1) {$\VV^{\tiny \sigma op}$};
\node (CT4) at (7,1) {$\otimes$};
\draw[->] (V1)--(W1) node[midway, left=-3pt]{$\Id$};
\draw[->] (V2)--(A2) node[midway, left=-3pt]{$\AA_\eta$};
\draw[->] (A2)--(W2) node[midway, left=-3pt]{$\overline{\AA}_\eta$};
\draw[->] (V3)--(A3) node[midway, left=-1pt]{$\AA_\eta$};
\draw[->] (A3)--(B3) node[midway, left=-3pt]{$\Id$};
\draw[->] (A'3)--(BT'3) node[midway, left = -3pt]{$coev_\VV$};
\draw[->] (B''3)--(C'3) node[midway, right=-3pt]{$\Id$};
\draw[->] (BT3)--(C3) node[midway, right=-3pt]{$ev_\VV$};
\draw[->] (C'3)--(W3) node[midway, left=-1pt]{$\overline{\AA}_\eta$};
\draw[->] (V4)--(AT4) node[midway, left=-3pt]{$coev_\VV$};
\draw[->] (A4)--(B4) node[midway, left=-3pt]{$\overline{\AA}_\eta$};
\draw[->] (B4)--(C4) node[midway, left=-3pt]{${\AA}_\eta$};
\draw[->] (A'4)--(C'4) node[midway, left=-3pt]{$\Id$};
\draw[->] (CT4)--(W4) node[midway, left=-3pt]{$ev_\VV$};
\draw[->] (V5)--(W5) node[midway, right=-3pt]{$\ZZ_\VV(S^1_{nb})$};
\node at (0.75,2.3) {$\MM_\eta$};
\node[scale = 1.4] at (0.75,1.9) {$\Rightarrow$};
\node at (2.5,2) {$\overset{snake}{\simeq}$};
\node at (5.5,2) {$\overset{sym.}{\simeq}$};
\node at (5.5,2) {$\overset{sym.}{\simeq}$};
\node at (8.75,2.2) {$\begin{array}{c} \Id_{coev_\VV}\\ \circ_1\\(\MM_T\otimes \Id_{\VV^{\tiny \sigma op}})\\ \circ_1\\ \Id_{ev_\VV} \end{array}$};
\node[scale = 1.4] at (8.75,1) {$\Rightarrow$};
\draw[blue!50, very thick, opacity = 0.5] (1.3,2.7) -- (1.3,1.3);
\draw[blue!50, very thick, opacity = 0.5] (3.9,3.4) ..controls (2.7,1) and (3.8,1)..(4,2);
\draw[blue!50, very thick, opacity = 0.5] (4.1,0.6) ..controls (5.3,3) and (4.2,3)..(4,2);
\draw[blue!50, very thick, opacity = 0.5] (6.6,1.5) ..controls (6.6,0) and (7.4,0)..(7.4,2) .. controls (7.4,3) and (6.5,3.5) .. (7,3.5).. controls (7.5,3.5) and (6.6,3).. (6.6,2.5);
\begin{scope}[xshift = 0.2cm, yshift = -0.4cm]
\draw[blue!50, very thick, opacity = 0.5] (9.7,2) ..controls (9.7,1) and (10.3,1)..(10.3,2) .. controls (10.3,3) and (9.5,3.5) .. (10,3.5).. controls (10.5,3.5) and (9.7,3).. (9.7,2);
\end{scope}
\end{tikzpicture}
\end{figure}
\\
Note that every bimodule above is induced by a functor as displayed here:\\
\begin{tikzpicture}[xscale = 1.3, yscale = 0.8]
\node (V1) at (0,1) {$\Vect_\FK$};
\node (V2) at (1.5,1) {$\VV \underset{\VV}{\boxtimes}\VV$};
\node (V3) at (4,1) {$(\VV \otimes \VV) \underset{\VV \otimes \VV^{\tiny \sigma op} \otimes \VV}{\boxtimes} (\VV \otimes \VV)$};
\node (V4) at (7,1) {$\underset{\Vect_\FK \otimes \VV^{\tiny \sigma op}}{\VV \ \boxtimes\ \VV}$};
\node (V5) at (10,1) {$\underset{\quad \VV \otimes \VV^{\tiny \sigma op}}{\VV \ \boxtimes\ \VV} $};
\node (W1) at (0,0) {$\FK$};
\node (W2) at (1.5,0) {$\unit \boxtimes \unit$};
\node (W3) at (4,0) {$(\unit \otimes \unit)\boxtimes (\unit \otimes \unit)$};
\node (W4) at (7,0) {$\unit \boxtimes \unit$};
\node (W5) at (10,0) {$\unit \boxtimes \unit$};
\draw[->] (V1)--(V2) node[midway, above=-2pt]{$\eta$};
\draw[->] (V2)--(V3) node[midway, above=-2pt]{$\sim$};
\draw[->] (V3)--(V4) node[midway, above=-2pt]{$\sim$};
\draw[->] (V4)--(V5) node[xscale = 0.75, yscale = 0.8,midway, above=-2pt]{\footnotesize $\Id_\VV \boxtimes (T \otimes \Id) \boxtimes \Id_\VV$};
\draw[|->] (W1)--(W2);
\draw[|->] (W2)--(W3);
\draw[|->] (W3)--(W4);
\draw[|->] (W4)--(W5);
\end{tikzpicture}\\
So $\RR_\VV(S^1_{nb})$ is induced by the monoidal functor given by inclusion of the unit in $\ZZ_\VV(S^1_{nb})$.

\section{Non-semisimple WRT relative to CY} \label{sectAppli}
We can now state the conjectures which are the main motivation for the study above. The main idea is that the Witten--Reshetikhin--Turaev theories and their non-semisimple variants can be obtained in a fully extended setting from a 3D theory relative to an invertible 4D anomaly. In particular, they are defined in a setting where the cobordism hypothesis applies, and can be rebuilt out of their value at the point. These would be a not necessarily semisimple modular tensor category for the invertible 4-TQFT and the 1-morphism induced by the inclusion of the unit for the relative 3-TQFT. As exposed above, in the non-semisimple case the unit inclusion is only partially dualizable, and induces a non-compact TQFT.

These conjectures follow ideas of Walker \cite{WalkerNotes}, Freed and Teleman \cite{FreedSlides} in the semisimple case, of Jordan and Safronov in the non-semisimple case. We do not know of a formal statement in the existing literature and propose one here.
\subsection{Bulk+Relative=Anomalous} \label{subsBRA}
Remember that the WRT theories, and their non-semisimple variants, are defined on a category of cobordisms equipped with some extra structure. They morally come from the data of a bounding higher manifold. 3-manifolds come equipped with an integer, which corresponds to the signature of the bounding 4-manifold, and surfaces come equipped with a Lagrangian in their first cohomology group, which corresponds to the data of the contractible curves in a bounding handlebody. In this setting, this extra structure is used to resolve an anomaly, and is due to Walker. We describe below how this kind of extra structure arises in the setting of relative field theories.
\begin{mydef}
The \emph{$(n-1)$-category of filled bordisms} $$\Bord_{n-1}^{filled}\subseteq \Bord_n^\to$$ is the sub-$(n-1)$-category of bordisms that map to the empty under the target functor $\Bord_n^\to \to \Bord_n$ and to $\Bord_{n-1}$ under the source functor. These are $k$-bordisms, $k \leq n-1$, equipped with a bounding $(k+1)$-bordism which we call the filling. We denote $$\Hollow:\Bord_{n-1}^{filled} \to \Bord_{n-1}$$ the functor that forgets the filling, namely the source functor. 

The \emph{$(n-1)$-category of non-compact filled bordisms} $$\Bord_{n-1}^{nc, filled}\subseteq \Bord_n^\to$$ is the sub-$(n-1)$-category of bordisms that map to the empty under the target functor and to $\Bord_{n-1}^{nc}$ under the source functor.
\end{mydef}
\begin{mydef}
An \emph{$n$-relative pair} $(\ZZ,\RR)$ is the data of:\begin{itemize}
    \item[]  an $n$-TQFT $\ZZ:\Bord_n\to \CC$
    \item[] an oplax-$\ZZ$-twisted-$(n-1)$-TQFT $\RR: \Bord_{n-1}\to \CC^\to$, namely an oplax transformation $\Triv \Rightarrow \ZZ\vert_{\Bord_{n-1}}$.
\end{itemize}
Such a pair is called a \emph{non-compact $n$-relative pair} if $\RR$ is a non-compact theory.
\end{mydef}
Given an $n$-relative pair $(\ZZ,\RR)$ one has two symmetric monoidal functors $\Bord_{n-1}^{filled} \to \CC^\to$. One is given by applying functoriality of $(-)^\to$ on $\ZZ$, namely applying $\ZZ$ to any diagram in $\Bord_n$ to get a diagram of the same shape in $\CC$. It has trivial target and gives an oplax transformation $$\ZZ^{\to \unit}: \ZZ\vert_{\Bord_{n-1}}\circ \Hollow \Rightarrow \Triv$$ between functors $\Bord_{n-1}^{filled} \to \CC$.

The other one is given by applying the relative field theory on the hollowed out bordism, it is an oplax transformation $$\RR \circ \Hollow: \Triv \Rightarrow\ZZ\vert_{\Bord_{n-1}}\circ \Hollow\ .$$
\begin{mydef}
The \emph{anomalous $(n-1)$-theory} $\AA$ induced by the $n$-relative pair $(\ZZ,\RR)$ is the composition $\ZZ^{\to \unit} \circ (\RR \circ \Hollow)$ of these two oplax transformations. It gives an oplax ransformation $\Triv \Rightarrow \Triv$ which by \cite[Theorem 7.4 and Remark 7.5]{JFS} is equivalent to a symmetric monoidal functor $$\AA: \Bord_{n-1}^{filled} \to (\Omega \CC)^{odd\ opp}\ ,$$ where $odd\ opp$ means we take opposite of $k$-morphisms for $k$ odd, and $\Omega\CC:= \End_\CC(\unit)$ is the delooping $(n-1)$-category. 

If $(\ZZ,\RR)$ is a non-compact $n$-relative pair, the same construction on the appropriate subcategories gives an anomalous theory $\AA: \Bord_{n-1}^{nc, filled} \to (\Omega\CC)^{odd\ opp}$.
\end{mydef}
For comparison with WRT theories, we will need to restrict to a once extended theory, namely look at endomorphisms of the trivial in $\Bord_{n-1}^{filled}$, and to check that the anomalous theory descends to the quotient where one only remembers signatures and Lagrangians out of the fillings. We will also move this odd opposite to the source category. 
\begin{mydef}
The \emph{bicategory of simply filled 3-2-1-cobordisms} $\cob_{321}^{filled}$ is the subcategory of $h_2(\Omega \Bord_3^{filled,\ odd\ opp})$ where circles can only be filled by disks, and surfaces by handlebodies. Taking the opposite orientation for 1- and 2-manifolds (which will have the effect of switching the source and target of a 3-bordism), one can identify this bicategory as:\\
$\cob_{321}^{filled} \simeq \left\{\begin{array}{ccl} \text{objects}&:& (\sqcup^n S^1,\sqcup^n \Dd^2:\sqcup^n S^1 \to \emptyset),\ n\in \Nn\\
\text{1-morphisms}&:& (\Sigma: \sqcup^{n_1} S^1 \to \sqcup^{n_2} S^1, H:  \emptyset \to (\sqcup^{n_1} \overline{\Dd^2}) \cup \Sigma \cup (\sqcup^{n_2} \Dd^2))\\
\text{2-morphisms}&:& (M: \Sigma_1 \to \Sigma_2, W: H_1 \cup M \cup \overline{H_2} \to \emptyset)
\end{array}\right.$\\
The analogous definition in the non-compact case $\cob_{321}^{nc, filled}\subseteq h_2(\Omega{\Bord_3^{nc, filled,\ odd\ opp}})$ will require 3-bordism to have non-empty \emph{incoming} boundary in every connected component, as source and targets of 3-manifolds are switched. To facilitate comparison with the existing literature, we also require that all surfaces have non-empty incoming boundary, although in our setting this is purely artificial. 
\end{mydef}
This bicategory is to be compared with:
\begin{mydef}
The bicategory $\widetilde{\cob}_{321}$ (resp. $\widetilde{\cob}_{321}^{nc}$) is the bicategory of circles, surfaces bordisms (resp. surface bordisms with non-empty incoming boundary) equipped with a Lagrangian subspace in their first homology group, and 3-bordisms (resp. 3-bordisms with non-empty incoming boundary) equipped with an integer. Composition is given by usual composition on the underlying bordisms, plus: 

$\scriptstyle \bullet\ $ taking the sum of the Lagrangian subspaces for composition of surfaces,

$\scriptstyle \bullet\ $ adding the integers plus some Maslov index for composition of 3-bordisms,

$\scriptstyle \bullet\ $ just adding the integers for composition of 3-bordisms in the direction of 1-morphisms.\\
See \cite[Section 3]{DeRenziDGGPR} and references therein for a precise definition. The bordisms there are decorated by objects of a ribbon category, and we are looking at the subcategory where every decoration is empty. The category $\widetilde{\cob}_{321}^{nc}$ corresponds to admissible bordisms there.
\end{mydef}
\begin{prop}
The assignment $$\pi_{321}:\left\{\begin{array}{ccc} \cob_{321}^{filled} &\to& \widetilde{\cob}_{321} \\ 
(\sqcup^n S^1,\sqcup^n \Dd^2) &\mapsto& \sqcup^n S^1\\
(\Sigma, H) & \mapsto & (\Sigma, \ker(i_*:H_1(\Sigma) \to H_1(H)))\\
(M,W) &\mapsto& (M, \sigma(W))
\end{array}\right.$$ is a symmetric monoidal functor. 
\end{prop}
\begin{myproof}
For composition of 1-morphisms we want to show that the kernel of a gluing is the sum of the kernels. One inclusion is immediate and the other one follows by dimensions since both are Lagrangians, see \cite[Propositions B.6.5 and B.6.6]{DeRenziPhD}. 

For composition of 2-morphisms we use Wall's theorem, see \cite[Theorem B.6.1]{DeRenziPhD} for a statement in our context.

For composition of 2-morphisms in the direction of 1-morphisms we use that filled surfaces only glue on disks, hence filled 3-manifolds on 3-balls, so the signature of the filling simply adds.
\end{myproof}
Similarly, one can restrict to non-compact cobordisms and get $$\pi_{321}^{nc}: \cob_{321}^{nc, filled}\to \widetilde{\cob}_{321}^{nc}\ .$$
If we restrict $\widetilde{\cob}_{321}$ to surfaces equipped with Lagrangians that are induced by some handlebody, these functors are essentially surjective, hence the name.

\subsection{Conjectures} \label{ssectConj}
We want to relate the Witten--Reshetikhin--Turaev theories and their non-semisimple variants to the ones induced by the cobordism hypothesis. We want to say that the anomalous theory induced the relative pair $(\ZZ_\VV,\RR_\VV)$ factors through $\widetilde\cob_{321}$ and recovers $\WRT$ and $\DGGPR$ theories. 

It has long been a folklore result that WRT theories extend to the circle \cite{WalkerOnWitten3mfldInv, GelcaTQFTwithCorners}, see also \cite{BalsamKirollovExtendedTV} for Turaev--Viro theories.
 Once-extended 3-TQFTs are classified in the preprint \cite[Theorem 3]{BDSPV}, and the following result can be obtained from \cite[Proposition 6.1]{BDSPV} (in our case the unit is simple). We give the statement of \cite[Theorem 1.1.1]{DeRenziPhD} restricted to trivially decorated bordisms.
\begin{myth}
For a semisimple modular tensor category $\mathscr{V}$ with a chosen square root of its global dimension, the Witten--Reshetikhin--Turaev TQFT extends to the circle as a symmetric monoidal functor $$\WRT_\mathscr{V}: \widetilde{\cob}_{321}\to \widehat Cat_\FK$$
where $\widehat Cat_\FK$ is the category of Cauchy-complete categories.
\end{myth}
Similarly, restricting the statement of \cite{DeRenziDGGPR} to trivially decorated bordisms:
\begin{myth}[Theorem 1.1 in \cite{DeRenziDGGPR}]
For a non-semisimple modular tensor category $\mathscr{V}$ with a chosen square root of its global dimension, the non-semisimple TQFT from \cite{DGGPR} extends to the circle as a symmetric monoidal functor $$\DGGPR_\mathscr{V}: \widetilde{\cob}_{321}^{nc}\to \widehat Cat_\FK$$
\end{myth}
On the other hand, using the Cobordism Hypothesis:
\begin{myth}[Brochier--Jordan--Safronov--Snyder]
For a semisimple or non-semisimple modular tensor category $\mathscr{V}$, its Ind-cocompletion $\VV \in \BrTens$ is 4-dualizable and induces under the Cobordism Hypothesis a 4-TQFT $\ZZ_\VV: \Bord_4^{fr}\to \BrTens$.
\end{myth}
The main result of this paper can be stated in this context.
\begin{myth}
For a semisimple modular tensor category $\mathscr{V}$, the arrow $\Aflat_\eta\in \BrTens^\to$ induced by the unit inclusion $\eta: \Vect_\FK \to \VV := Ind(\mathscr{V})$ is 3-dualizable and induces under the Cobordism Hypothesis a framed oplax-$\ZZ_\VV$-twisted 3-TQFT $$\RR_\VV: \Bord_3^{fr}\to \BrTens^\to\ .$$
For a non-semisimple modular tensor category $\mathscr{V}$, $\Aflat_\eta$ is not 3-dualizable but is non-compact-3-dualizable and induces under the non-compact Cobordism Hypothesis a framed non-compact oplax-$\ZZ_\VV$-twisted 3-TQFT $$\RR_\VV: \Bord_3^{fr, nc}\to \BrTens^\to\ .$$
\end{myth}
\begin{myproof}
If $\Vv$ is semisimple, $\VV=Ind(\Vv)=Free(\Vv)$ and Theorem \ref{thmCrit3oplax} applies. If $\Vv$ is not semisimple, the unit is not projective in $\Vv$, nor in $\VV=Ind(\Vv)$, so $\Aflat_\eta$ is not 3-dualizable. But $\VV$ is cp-rigid and Theorem \ref{thmCritnc3d} applies.
\end{myproof}
To compare the two sides, we need all theories to be oriented. We assume the following: 
\begin{conjecture}\label{conjSO3} Let $\Vv$ be a ribbon tensor category and $\VV:=Ind(\Vv)$, then: \\
  The ribbon structure of $\mathscr{V}$ induces an $SO(3)$-homotopy-fixed-point structure on $\VV$.\\
  The ribbon structure of $\eta$ induces an $SO(3)$-homotopy-fixed-point structure on $\Aflat_\eta$.
\end{conjecture}
The first statement is expected by experts. The second one follows \cite[Example 4.3.23]{LurieCob}. Note that in the second statement we really mean an $SO(3)$-homotopy-fixed-point structure compatible with the one on $\VV$, as in Remark \ref{rmkOrientedTwistedCH}.
\begin{rmk} The fact that the anomalous theory $\AA_\VV$ would factor through $\widetilde{\cob}_{321}$ is not too surprising. As was pointed to me by Pavel Safronov, we know from \cite{BJSS} that $\VV$ is not only 4-dualizable, but invertible, and hence 5-dualizable. But $ \BrTens $ has no non-trivial 5-morphisms, and hence the 5-theory induced by $\VV$ is trivial on 5-bordisms. This means that $\ZZ_\VV$ should give the same value on cobordant 4-manifolds. If this story can be made oriented, it means it depends only on the signature of 4-manifolds.

It was observed by Walker \cite[Chapter 9]{WalkerNotes} in the semisimple case that there is a scalar choice of ways to extend $\ZZ_\VV$ from $\Bord_3^{or}$ to $\Bord_4^{or}$, namely $\ZZ_\VV(B^4)$, and that exactly two of these scalars yield theories which are cobordant-invariant on 4-manifolds. He observes that these scalars are exactly the two square roots of the global dimension among which one has to choose when defining WRT theories.
This motivates the following conjecture. In the non-semisimple case, it is supported by the fact that the constructions of the (3+1)-TQFTs of \cite{CGHP} need exactly the choice of a modified trace.
\end{rmk}
\begin{conjecture}\label{conjSO45} Let $\Vv$ be a modular tensor category and $\VV:=Ind(\Vv)$, then: \\
  A choice of modified trace on $\mathscr{V}$ induces an $SO(4)$-homotopy-fixed-point structure on $\VV$.\\
A modified trace induces an $SO(5)$-homotopy-fixed-point structure on $\VV$ if and only if the global dimension $\mathscr S_{\Vv, \mt}(S^4)=1$ with this choice of modified trace in the construction of \cite{CGHP}.
\end{conjecture}
In particular, we conjecture that every modular tensor category has an $SO(5)$-homotopy-fixed-point structure. Indeed let $\mathscr V$ be a modular tensor category and choose $\mt$ a non-degenerate modified trace on $\mathscr V$, which exists and is unique up to scalar by \cite[Corollary 5.6]{GKPmtrace}. Choose a square root $\mathscr D_\mt$ of its global dimension $d(\Vv)_\mt := \mathscr S_{\mathscr V,\mt}(S^4) = \Delta_+\Delta_-$ as defined in \cite{CGHP} (and denoted $\zeta$ there). Then the modified traces $\pm\mathscr D_\mt^{-1} \mt$ are the only two modified traces satisfying $\mathscr S_{\mathscr V, \pm\mathscr D_\mt^{-1} \mt}(S^4)=1$ by \cite[Proposition 5.7]{CGHP}.
\begin{rmk}
Let us try to give a conceptual reason for why $SO(5)$-structures correspond to square roots of the global dimension. As $\VV$ is an invertible object, the oriented theory $\ZZ_\VV$ is an invertible 4-TQFT and these are known to give invariants which only depend on the signature and Euler characteristic of 4-manifolds. Two closed 4-manifold are cobordant if and only if they have the same signature, so to get a cobordant-invariant theory we need to kill the dependence on the Euler characteristic. Changing the choice of modified traces by a scalar $\kappa$ alters this dependence by a factor $\kappa^{\sigma(W)}$ \cite[Proposition 5.7]{CGHP}\footnote{We assume Conjecture \ref{conjCGHP} here.}, and a well-chosen scalar $\kappa$ will kill it. We need to ask that for any closed 4-manifold $W$ with signature zero, $\ZZ_\VV(W)=1$. However there is no closed 4-manifold with signature 0 and Euler characteristic 1, they always have same parity. It is sufficient to ask that $\ZZ_\VV(S^4)=1$. The 4-sphere has Euler characteristic 2, hence there are exactly two solutions for $\kappa$, the two square roots of the global dimension.
\end{rmk}
\begin{corr}[of conjectures]
Both $\ZZ_\VV$ and $\RR_\VV$ give oriented TQFTs by the oriented cobordism hypothesis.  
\end{corr}
We now assume that this corollary is true, that the choice of square root of the global dimension has been made, and that $\ZZ_\VV$ and $\RR_\VV$ are oriented.

In the semisimple case, the relative pair $$(\ZZ_\VV:\Bord_4\to \BrTens,\RR_\VV:\Bord_3\to \BrTens^\to)$$ induces an anomalous theory $$\AA_\VV: \Bord_3^{filled,\ odd\ opp} \to \Tens := \Omega{\BrTens}\ .$$ Its restriction on $\cob^{filled}_{321}$ gives a 2-functor $$\AA_\VV^{321}: \cob^{filled}_{321} \to \Omega{\Tens}\simeq \Pr\ .$$
\begin{conjecture}\label{conjSS}
For a semisimple modular tensor category $\mathscr{V}$, the anomalous theory induced by $(\ZZ_\VV,\RR_\VV)$ recovers the Witten--Reshetikhin--Turaev theory. Namely: \center
\begin{tikzpicture}[xscale = 4, yscale = 2, baseline = 5pt]
\node (Cobfilled) at (0,1) {$\cob^{filled}_{321}$};
\node (Cobtilde) at (0,0) {$\widetilde{\cob}_{321}$};
\node (Cat) at (1,0) {$\widehat Cat_\FK$};
\node (Pr) at (1,1) {$\Pr$};
\draw[->] (Cobfilled) -- (Cobtilde) node[midway, left]{$\pi_{321}$};
\draw[->] (Cobfilled) -- (Pr) node[midway, above]{$\AA_\VV^{321}$};
\draw[->] (Cobtilde) -- (Cat) node[midway, above]{$\WRT_\mathscr{V}$};
\draw[right hook->] (Cat) -- (Pr) node[midway, right]{$\Free$};
\end{tikzpicture} \quad commutes up to isomorphism.
\end{conjecture}
In the non-semisimple case, the relative pair $$(\ZZ_\VV:\Bord_4\to \BrTens,\RR_\VV:\Bord_3^{nc}\to \BrTens^\to)$$ induces an anomalous theory $$\AA_\VV: \Bord_3^{nc, filled,\ odd\ opp} \to \Tens := \Omega{\BrTens}\ .$$ Its restriction on $\cob^{nc,filled}_{321}$ gives a 2-functor $$\AA_\VV^{321}: \cob^{nc,filled}_{321} \to \Omega{\Tens}\simeq \Pr\ .$$
\begin{conjecture}\label{conjNonSS}
For a non-semisimple modular tensor category $\mathscr{V}$, the non-compact anomalous theory induced by $(\ZZ_\VV,\RR_\VV)$ recovers the De Renzi--Gainutdinov--Geer--Patureau-Mirand--Runkel theory. Namely: \center
\begin{tikzpicture}[xscale = 4, yscale = 2, baseline = 5pt]
\node (Cobfilled) at (0,1) {$\cob^{nc, filled}_{321}$};
\node (Cobtilde) at (0,0) {$\widetilde{\cob}_{321}^{nc}$};
\node (Cat) at (1,0) {$\widehat Cat_\FK$};
\node (Pr) at (1,1) {$\Pr$};
\draw[->] (Cobfilled) -- (Cobtilde) node[midway, left]{$\pi_{321}^{nc}$};
\draw[->] (Cobfilled) -- (Pr) node[midway, above]{$\AA_\VV^{321}$};
\draw[->] (Cobtilde) -- (Cat) node[midway, above]{$\DGGPR_\mathscr{V}$};
\draw[right hook->] (Cat) -- (Pr) node[midway, right]{$\Free$};
\end{tikzpicture} \quad commutes up to isomorphism.
\end{conjecture}
We know how to check these conjectures on the circle. We have $\WRT_\mathscr{V}(S^1) = \Vv$ whose free cocompletion is equivalent to $\VV$ because $\Vv$ is semisimple. Similarly, $\DGGPR_\mathscr{V}(S^1)=\Proj(\mathscr{V})$ whose free cocompletion is equivalent to $\VV$.
 On the other side, we know that in dimension two $\ZZ_\VV$ coincides with factorization homology, and we computed $\RR_\VV(S^1)$ in Section \ref{subsCircle}. So: $$\AA_\VV^{321}(S^1,\Dd^2) = \RR_\VV(S^1) \underset{\ZZ_\VV(S^1)}{\boxtimes}\ZZ_\VV(\Dd^2) \simeq {}_{\Vect_\FK}\ZZ_\VV(S^1) \underset{\ZZ_\VV(S^1)}{\boxtimes}\VV_{\Vect_\FK} \simeq {}_{\Vect_\FK}\VV_{\Vect_\FK}$$

Computing the values of the theories induced by the Cobordism Hypothesis on higher dimensional bordisms comes down to computing some adjoints in $\BrTens$ and compose them in various ways. This will be carried out in future work.
\begin{corr}[of conjectures]
Both $\WRT_\Vv$ and $\DGGPR_\Vv$ extend to $S^0$.
\end{corr}
\begin{myproof}
Indeed, the anomalous theory $\AA_\VV$ is really defined as a functor between the 3-categories $\Bord_3^{filled} \to\Tens$ (resp. $\Bord_3^{nc, filled} \to\Tens$ in the non-semisimple case). The two points $S^0$ are bordant, by a cap, and therefore give an object $(S^0,\cap) \in \Bord_3^{filled}$ (resp. $\Bord_3^{nc, filled}$). 

It is easy to compute the value of the anomalous theory on this object, namely $\AA_\VV(S^0,\cap) = \RR_\VV(S^0) \circ \ZZ_\VV(\cap) = (\AA_\eta \boxtimes (\overline{\AA}_\eta)^*) \underset{\VV \boxtimes \VV^{\sigma op}}{\boxtimes} \VV \simeq \VV$ seen as a $\Vect_\FK$-$\Vect_\FK$-central algebra.
\end{myproof}
\begin{rmk}
This corollary is to be compared with results of \cite{DSPS} which shows that $\WRT_\Vv$ extends to the point if and only if $\Vv \simeq Z(\mathscr{C})$ is a Drinfeld center, in which case the point is mapped to $\mathscr{C}$. In the modular case, the Drinfeld center $Z(\mathscr{C})$ is isomorphic to $\mathscr{C} \otimes \mathscr{C}^{\sigma op}$, and the two descriptions agree on $S^0$. Therefore it appears that $\WRT_\Vv$ always extends to $S^0$, and extends to the point if and only if one can find a ``square root" for its value on $S^0$. This is also related to ongoing work of Freed, Teleman and Scheimbauer.

Note however that the statement above is a bit informal, because it is really $\Free\circ \WRT_\Vv \circ \pi_{321}$ that extends to $S^0$, so WRT indeed but with different source and target. In particular, the results of \cite{DSPS} do not apply directly in this context.
\end{rmk}

\newpage \footnotesize
\bibliography{mybib}{}
\bibliographystyle{alpha}

\end{document}